\documentclass{elsart}
\usepackage{amsmath,amssymb,verbatim,alltt,amsfonts,amsfonts,array,texdraw}
\usepackage{latexsym, xy,amscd,url,epic}
\usepackage[english]{babel}
\journal{}

 \DeclareMathOperator{\Co}{Coker}
 \DeclareMathOperator{\im}{Im}
 
 \newcommand{\To}{\longrightarrow}
 \newcommand{\ra}{\rightarrow}
 \newcommand{\km}{\mathcal{M}}
 \newcommand{\ke}{\mathcal{E}}
 \newcommand{\kf}{\mathcal{F}}
 \newcommand{\kb}{\mathcal{B}}
 \newcommand{\kl}{\mathcal{L}}
 \newcommand{\kc}{\mathcal{C}}

 \newcommand{\Z}{\mathbb{Z}}

\begin{document}

\begin{frontmatter}
\theoremstyle{definition}

\title{\textbf{Maximal Cohen-Macaulay Modules over the Affine Cone of the Simple Node}} 

\author{Corina Baciu} 

\ead{baciu@mathematik.uni-mainz.de}

\address{Johannes Gutenberg-Universität Mainz,
Fachbereich Physik, Mathematik und Informatik, Staudingerweg 9, D-55099 Mainz}

\begin{abstract}  A concrete description of all graded maximal Cohen--Macaulay modules of rank one and two 
over the non-isolated singularities of type $y_1^3+y_1^2y_3-y_2^2y_3$ is given. For this purpose we construct an alghoritm that provides extensions of MCM modules over an arbitrary hypersurface.
\end{abstract}

\begin{keyword} Hypersurface ring, Maximal Cohen--Macaulay modules, Non-isolated singularity \\$2000$ $Mathematics$ $Subject$ $Classification$: 13C14, 13H10, 14H60, 14H45, 16W50, 32S25
\end{keyword}

\end{frontmatter}

\section*{Introduction}

Over a (graded) hypersurface ring, the (graded) maximal Cohen Macaulay modules (shortly MCM modules) can be described 
by (graded) matrix factorizations of the polynomial that defines the hypersurface. (see \cite{Ei}) 

Let $S$ be a polynomial ring over a field $k$ and $f\in S$ an irreducible homogeneous polynomial  of degree $d$. 
Consider the hypersurface ring $R=S/f$  and $M$ a graded MCM module over it.
D. Eisenbud proved that $M$ has an infinite, graded, 2-periodic resolution  over $R$, of the form
  \[
...\stackrel{\varphi}{\longrightarrow}\overset{j=n}{\underset{j=1}{\oplus}}R(\alpha_{j}-d)
\stackrel{\psi}{\longrightarrow}\overset{j=n}{\underset{j=1}{\oplus}}R(\beta_{j})
\stackrel{\varphi}{\longrightarrow}\overset{j=n}{\underset{j=1}{\oplus}}R(\alpha_{j})\longrightarrow M
\longrightarrow 0,   \]
where the  maps $\varphi$ and $\psi$ are the multiplications by some square  matrices $A$, $A'$ with homogeneous entries and fulfill 
$\varphi\psi=\psi\varphi=fId$. We say that the pair $(A,A')$ form a graded $matrix$ $factorization$ of the polynomial $f$.\\ 
Notice that the matrix factorization alone determine the module only up to shiftings. In order to obtain also some information on the degree, we have to know the coefficients $\alpha_{j}$ and $\beta_{j}$. The rank of the module $\textrm{Coker}\,\varphi$ is precisely the integer $r$ such that $\textrm{det}\,A=f^r$. It follows immediately that the minimal number of generators of a MCM $R$--module of rank $r$ is smaller equal $dr$.
                   
 If $k$ is an algebraically closed field, the graded MCM $R$--modules can be described also geometrically. 
They correspond, by sheafification, to the aCM ($arithmetically$ $Cohen$--$Macaulay$) sheaves on the projective cone
$Y=\textrm{Proj}\, S$ (Grothendieck, see for ex. \cite{Groth}). If $Y$ is a smooth curve, the aCM sheaves are exactly the vector bundles.

In \cite{LPP}, using the classification of the vector bundles over a smooth elliptic curve realized by Atiyah (see \cite{At}), the authors produced an alghoritm to construct matrix factorizations of all indecomposable graded MCM modules over $R=k[y_1,y_2,y_3]/\langle y_1^3+y_2^3+y_3^3\rangle$. Essential is that a smooth elliptic curve has a tame category of vector bundles.\\
In \cite{DG},  G.\,-\,M.\,Greuel and Yu.\,Drozd proved that the nodal curve has also a tame category of vector bundles, even a tame category of coherent sheaves. 

The aim of this paper is to describe explicitly all graded MCM modules of low ranks (1 and 2) over the affine cone of the nodal curve.\\The matrix factorizations of all graded, rank one MCM modules are constructed in section 2. 
The rank two non--locally free aCM sheaves are described as extensions of rank one aCM sheaves (see section 4). In fact, with the same method, inductively, one can construct  indecomposable non--locally free aCM sheaves of any rank.\\
 The main tool is the theorem $\ref{mfext}$ that gives us a way to construct a matrix factorization of a module $E$, that fits in a graded extension of type 
\[
0\ra L \ra E \ra F \ra 0
\] 
where $L$ and $F$ are two graded MCM modules with known matrix factorizations. 
The  computations are made with the help of the computer algebra system \textsc{Singular}.(see \cite{GPS})\\
We obtain 4 families parametrized by the regular points of the curve and 22 countable families.
(Remark: It is known that a graded hypersurface has countable CM--representation type iff it is isomorphic to $A_\infty$; 
see \cite{Trau})\\
As a direct application of the theorem $\ref{mfext}$ we prove that, over the affine cone of curves of arithmetic genus 1, the matrix factorizations  give information on the stability of the sheafification of graded MCM modules.\\
In the third section, we use the classification of the vector bundles on the simple node (see \cite{DG}) in order to give  matrix factorizations of the locally free rank two graded MCM modules. One can construct them also using the extensions of rank one locally free MCM modules, but knowing the form of the corresponding vector bundle on the simple node, one can get some extra geometrical properties. For example, one can find the matrix factorizations corresponding to all stable rank two vector bundles. (see \cite{Chisinau} and \cite{thesis}) 
Other interesting ways of computing matrix factorizations can be found in \cite{BEPP}.

The author thanks Prof. Dr. V. Vuletescu  and  Dr. I. Bourban for helpful discussions on Section 3, respectively  on the classification of the vector bundles on the nodal curve. The author is grateful also to Prof. G. Pfister for helpful remarks on the Section 2.\\
The research was supported by the DFG-Schwerpunkt "Globale Methoden in der komplexen Geometrie".


\vspace{-0,1cm}

\section{Extensions of MCM modules over hypersurface rings}

In the following we show how to construct extensions of two MCM modules with known matrix factorizations.
This method will be used in the last section of this paper for the classification of all non--locally free rank two MCM modules over the affine cone of the simple node.\\ 

\begin{thm}{\label{mfext}}
Let $S=k[x_1,...,x_n]$ where $k$ is a field and let $R=S/f$  be a hypersurface ring defined by an irreducible homogeneous polynomial $f$.
Consider $L,F$ two graded MCM $R$--modules with the matrix factorizations $(A,A')$, respectively $(B,B')$
and  the graded extension $0\ra L \stackrel{\alpha}{\longrightarrow} E
\stackrel{\beta}{\longrightarrow} F \ra 0.$\\
 Then $E$ is a graded MCM $R$--module and has a
matrix factorization $(M,M')$, $M$ of the type $M=\left(\begin{matrix} A & D \\0 & B \end{matrix}\right)$,
$D$ a matrix with homogeneous entries in $S$ such that $A'\cdot D\cdot B'=0$ in $R$.
\end{thm}

\begin{pf}
Denote $s=\mu(L), t=\mu(F), r_1=$rank$(L), r_2=$rank$(F)$.\\
Consider the following graded diagram:
\[
\begin{CD}
0 \To \overset{j=s}{\underset{j=1}{\oplus}}R(\alpha'_{j})@>i>>\overset{j=s+t}{\underset{j=1}{\oplus}}R(\alpha'_{j})@>\pi>>\overset{j=s+t}{\underset{j=s+1}{\oplus}}R(\alpha'_{j}) @>>>0\\
\hspace{1cm}@VV{A}V     @VVV  @VV{B}V\\
0 \To  \overset{j=s}{\underset{j=1}{\oplus}}R(\alpha_{j})@>i>>\overset{j=s+t}{\underset{j=1}{\oplus}}R(\alpha_{j})@>\pi>>\overset{j=s+t}{\underset{j=s+1}{\oplus}}R(\alpha_{j}) @>>>0\\
\hspace{1cm}@VV{p_A}V      @VV{\delta}V      @VV{p_B}V\\
0 \To \Co A @>\alpha>> E @>\beta>> \Co B @>>>0\\
\end{CD}
\]
where $\pi$ is the projection on the last $t$ variables and $i$ is the natural inclusion.
Let  $\varphi_B : \overset{j=s+t}{\underset{j=s+1}{\oplus}}R(\alpha_{j}) \to E$ be a graded map such that $\beta \circ \varphi_B =p_B$ and consider $\varphi_A = \alpha \circ p_A$.
Define the graded map $\delta$ as the sum of $\varphi_A$ and $ \varphi_B$. Then, $\delta$ makes the above diagram commutative.
 Using Snake Lemma, we get the graded exact sequence
\[
0\To \im A \stackrel{i}{\To} \textrm{Ker}\, \delta \stackrel{\pi}{\To} \im B \To 0
\]
and the surjectivity of $\delta$.
Therefore, $\textrm{Ker}\, \delta$ is a graded MCM $R$--module of rank $s+t-r_1-r_2$ such that $E\simeq \overset{j=s+t}{\underset{j=1}{\oplus}}R(\alpha_{j})/\textrm{Ker}\,
\delta$.\\
So, there exists $(M_1,M'_1)$ a graded matrix factorization of $f$ such that $E\simeq \Co M'_1$ and such that the following graded sequence
is exact$$
0\To \im A \stackrel{i}{\To} \im M'_1 \stackrel{\pi}{\To} \im B \To 0.$$

Since $\im \left(\begin{smallmatrix} A\\0\end{smallmatrix}\right)\subset \im M'_1$, there exist two invertible matrices $U$ and $V$ such that
$UM'_1V=\left(\begin{smallmatrix} A & D_1 \\0 & B_1\end{smallmatrix}\right)$. Denote $M''=UM'_1V$. Then $\Co M''$ is a graded MCM module isomorphic to
$E$. Thus, there exists the matrix $B'_1$  with homogeneous entries such that $(B_1,B'_1)$ is a matrix factorization of $f$
and $A'D_1B'_1=0$ in $R$.\\
We have therefore a graded commutative diagram:
\[
\begin{CD}
0 \To \im A @>i>>\im M''@>\pi>>\im B_1 @>>>0\\
\hspace{1cm}@VV{id}V     @VV{\wr}V  \\
0 \To \im A @>i>>\im M'_1@>\pi>>\im B @>>>0\\
\end{CD}
\]
and an inclusion $\im B_1\hookrightarrow \im B$ of two graded MCM modules of the same rank. This means that
they are isomorphic, and therefore, there exist two invertible matrices $U_1$ and $V_1$ such that $B=U_1B_1V_1$.
Denote $D=D_1V_1$ and $M=\left(\begin{smallmatrix} A & D \\0 & B \end{smallmatrix}\right)$. Then $A'DB'=0$ in $R$ and $\Co M \simeq \Co M'' \simeq E$.
 \end{pf}

So, if we have two graded matrix factorizations $(A,A')$, respectively $(B,B')$ (together with the corresponding graded maps) of two MCM $R$--modules
and we want an extension from $\text{Ext}^1(\Co B, \Co A)$ we do the following:\\
1) We decide the degree of the entries of a matrix $D$ such that $\left(\begin{smallmatrix} A & D \\0 & B \end{smallmatrix}\right)$ defines a graded map.
With the notations from the previous proof, the entry $(i,j)$ of $D$ should have the degree $\alpha_i-\alpha_{s+j}'$. Put a zero entry where the degree is negative. If all entries have negative degree, then we have only the zero extension.\\
2)We make some linear transformations in order to simplify the entries of $D$.\\
3)We impose the condition $A'DB'=0$ in $R$, that means all entries of $D$ reduce to zero modulo $f$. (often is necessary the computer)\\
4)We check if the obtained extension is non--zero, that means, do not exist matrices $U$ and $V$ such that $D=AU+VB$. \\
Remark: Sometimes it is worthy to ``inspire'' from the result at third step before to end the second one, since some transformations are more useful then another.

Concrete examples can be found in the last section of this article.\\
The reverse of theorem $\ref{mfext}$ is also true.
   \begin{prop}\label{mfext1}
Let $S$, $R$, $f$, $L$ and $F$  as in $\ref{mfext}$. 
   Let $D$ be a matrix with homogeneous entries in $S$ such that $A'\cdot D\cdot B'=0$ in $R$.\\
   Then there exists a graded MCM module $E$ with a matrix factorization $(M,M')$, $M=\left(\begin{matrix} A & D \\0 & B \end{matrix}\right)$ and an extension $   0\ra L \ra E \ra F \ra 0.$
   \end{prop}
   \begin{pf}
   Denote $s=\mu(L), t=\mu(F), r_1=$rank$(L), r_2=$rank$(F)$.\\ We know that $AA'=A'A=Id_s$ and $BB'=B'B=Id_t$.\\
   The condition $A'\cdot D\cdot B'=0$ modulo (f) means that there exists the matrix $C$ such that $DB'=AC$ and $A'D=CB$.
   Then $\left(\begin{smallmatrix} A & D \\0 & B \end{smallmatrix}\right)\cdot \left(\begin{smallmatrix} A' & -C \\0 & B'
   \end{smallmatrix}\right)=fId_{s+t}$ and
     $\left(\begin{smallmatrix} A' & -C \\0 & B' \end{smallmatrix}\right)\cdot\left(\begin{smallmatrix} A & D \\0 & B
   \end{smallmatrix}\right)=fId_{s+t}$.
  So, $\Co M$, for $M=\left(\begin{smallmatrix} A & D \\0 & B \end{smallmatrix}\right)$ is a MCM $R$--module.
   From the commutative diagram:
   \[
   \begin{CD}
   0 \longrightarrow R^s @>i>>R^{s+t} @>\pi>>R^t @>>>0\\
   \hspace{1cm}             @VV{A}V      @VV{M}V      @VV{B}V\\
   0  \longrightarrow R^s @>i>>R^{s+t} @>\pi>>R^t @>>>0\\
   \end{CD}
   \]
   with $\pi$ the projection on the last $t$ components, using the Snake--Lemma, we get $0\ra L \ra E \ra F \ra 0.$
   \end{pf}

As a direct corollary of Theorem $\ref{mfext}$, we give a description of the modules with stable sheafification on a projective curve with arithmetic genus 1.

Let $Y \subset \mathbb P^2$ be a rational curve. For a locally free sheaf on $Y$, $\mathcal{E}$, we define
 $\text{deg}\ \mathcal{E}$ to be $\text{deg}\ \mathcal{E}=\chi(\mathcal{E})+(p_a(Y)-1)\text{rank}\ (\mathcal{E})$.\\
The vector bundle $\mathcal{E}$ is called $stable$ if for any torsion free quotient
$\mathcal{F}$ of $\mathcal{E}$,
\begin{displaymath}
  \frac{\text{deg}\ \mathcal{E}}{\text{rank}\ \mathcal{E}}<\frac{\text{deg}\ \mathcal{F}}{\text{rank}\ \mathcal{F}}.
\end{displaymath}
A graded MCM module over the affine cone of $Y$ is called $stable$ if the its sheafification is a stable vector bundle on $Y$.

Denote $S=k[y_1,y_2,y_3]$ and let $F$ be the homogeneous polynomial defining the curve $Y$.\\
Let $M$ be a graded, indecomposable, locally free MCM module over the affine cone of $Y$ and $\mu$ the minimal number of generators of M.\\ 
Consider $(A,A')$ a graded matrix factorization of $F$ together with the map 
$\overset{i=\mu}{\underset{i=1}{\oplus}}R(\alpha_{i}) \stackrel{A}{\To} \overset{i=\mu}{\underset{i=1}{\oplus}}R(\beta_{i})$, corresponding to $M$. \\
We denote with $\mathcal{M}^{\alpha,\beta}_{\mu \times \mu}(S)$ the set of $\mu \times \mu$ matrices with a homogeneous entry of degree $\beta_i-\alpha_j$ on the position $(i,j)$, for all $i,j=1, \dots, \mu$.

Set $\kl_A $ the vector space
$$ \kl_A=\{D\in\mathcal{M}^{\alpha,\beta}_{\mu \times \mu}(S) | A'DA'=0\  \textrm{modulo} (F)\}$$
and define the following equivalence on it: two matrices $D$ and $D'$ from  $\kl_A $
are equivalent ($D \sim D'$) iff there exist two quadratic matrices $U$ and $V$ such that $D-D'=UA-AV$. Denote $\mathcal{S}_A=\kl_A/\sim$.\\

\begin{rem}
If $(B,B')$ is another matrix factorization of $M$, the vector spaces $\mathcal{S}_A$ and
$\mathcal{S}_B$ are isomorphic.
\end{rem}
Indeed, if $U$ and $V$ are the invertible matrices such that $B=UAV$, we construct the vector spaces isomorphism:
$\theta :\mathcal{S}_B \To \mathcal{S}_A, \theta(D)=U^{-1}DV^{-1}.$\\

\begin{thm}
Let $Y\subset \mathbb P^2$ be a projective curve  with arithmetic genus 1 and let $M$ be a graded, indecomposable,
locally free MCM module over the affine cone of $Y$.
The following statements are equivalent:
\begin{enumerate}
\item  $M$ is a stable module;
\item  $\textrm{dim}\ \mathcal{S}_A=1$ for $(A,A')$ a matrix factorization of $M$;
\item $\textrm{dim}\ \mathcal{S}_A=1$ for all $(A,A')$ matrix factorizations of $M$.
\end{enumerate}
\end{thm}

\begin{pf}
A  vector bundle $\ke$ on a curve with  arithmetic genus 1 is stable if and only if is simple, that means $\text{dim(Ext}^1(\ke,\ke))$=1 (see a proof in \cite{igor}).
The theorem $\ref{mfext}$ implies that $\text{dim(Ext}^1(\ke,\ke))=\textrm{dim}\ \mathcal{S}_A$ for $(A,A')$ a matrix factorization of $M$, so the first two statements are equivalent. The previous remark implies the equivalence of the last two statements.  
\end{pf}


\section{  Rank one, graded, MCM modules }

We know that the minimal number of generators of a rank one graded MCM $R$--modules is smaller equal to 3, the degree of $f$. 
Therefore, they are two or three minimally generated.

The line bundles on the simple node $Y=\textrm{Proj}R$ corresponding to the locally free modules are described in \cite{DG} as $\mathcal B(d,1,\lambda)$ with $\lambda \in k^*$ ($\lambda$ run over all regular points of
the curve $Y$) and $d$ the degree of the bundle. \\
The tensor product of two line bundles is given by:
\[ \mathcal B(d,1,\lambda)\otimes
\mathcal B(d',1,\lambda ')=\mathcal B(d+d',1,\lambda\cdot \lambda').\]

\subsection{Two--generated graded MCM R--modules}

   Let $s=(0:0:1)$ be the unique singular point of the curve $V(f)\subset \mathbb P^{2}_k$ and denote
$V(f)_{\mathrm{reg}}=V(f)\backslash \{s\}$.\\
Then $V(f)_{\mathrm{reg}}=\{(\lambda_1:\lambda_2:1),\lambda_1^3+\lambda_1^2-\lambda_2^2=0,\lambda_1\ne 0 \} \cup
\{(0:1:0) \}$.

  For any $\lambda=(\lambda_1:\lambda_2:1)$ in $V(f)$ denote:
\[\varphi_\lambda=\left(\begin{matrix}
  y_1-\lambda_1y_3 & \quad y_2y_3+\lambda_2y_3^2
\\y_2-\lambda_2y_3 & \quad y_1^2+(\lambda_1+1)y_1y_3+(\lambda_1^2+\lambda_1)y_3^2
\end{matrix}\right),
\]\[
\psi_\lambda=\left(\begin{matrix}
   y_1^2+(\lambda_1+1)y_1y_3+(\lambda_1^2+\lambda_1)y_3^2 & \quad -(y_2y_3+\lambda_2y_3^2)
\\      -(y_2-\lambda_2y_3)                               & \quad     y_1-\lambda_1y_3
                    \end{matrix}\right).
\]
 If $\lambda=(0:1:0)$ let be:
\[
\varphi_\lambda= \left( \begin{matrix}
   y_1+y_3 & y_2^2
\\ y_3 & y_1^2              \end {matrix}\right),\hspace{1cm}
\psi_\lambda= \left( \begin{matrix}
   y_1^2 & -y_2^2
\\ -y_3 & y_1+y_3            \end {matrix}\right).
\]
For any $\lambda \in V(f)_\mathrm{reg}$, we consider also the following graded maps defined by the above matrices:\\
$\psi_\lambda : R(-2)^2 \longrightarrow R \oplus R(-1)$ and  $\varphi_\lambda : R(-2) \oplus R(-3) \longrightarrow R(-1)^2$.\\\\
Define  $\mathcal M_{-1}=\{\,\mathrm{Coker}\,\varphi_\lambda\,|\,\lambda\in V(f)_\mathrm{reg}\,\},
\mathcal M_1=\{\,\mathrm{Coker}\,\psi_\lambda\,|\,\lambda\in V(f)_\mathrm{reg}\,\}$
and $\underline{\mathcal M}=\{\,\mathrm{Coker}\,\varphi_s, \mathrm{Coker}\,\psi_s\}$.\\

\begin{thm} \label{2gen}
(1) For all $\lambda \in V(f)$, $(\varphi_\lambda,\psi_\lambda)$ is a matrix factorization of $f$;  \\
(2) Every two-generated non-free graded MCM $R$-module is isomorphic, up to shifting, with one of the
modules from $\mathcal M_{-1} \cup \mathcal M_1 \cup \underline{\mathcal M}$; \\
(3) Every two different $R$-modules from
 $\mathcal M_{-1} \cup \mathcal M_1 \cup \underline{\mathcal M}$ are not isomorphic;\\
(4) All the modules from
$\mathcal M_{-1}\cup\mathcal M_1 \cup \underline{\mathcal M}$ have rank 1; \\
(5) The modules from $\mathcal M_1$ are the syzygies and also the duals of the modules from
$\mathcal M_{-1}$.
  \end{thm}
 \begin{pf}
(1) Since $\varphi_\lambda\psi_\lambda=\psi_\lambda\varphi_\lambda=f\cdot \mathbf 1_2$
for any $\lambda \in V(f)$, the first statement is true.

(2) Let $M$ be a two-generated non-free graded MCM $R$-module and consider $(\varphi,\psi)$ a graded reduced matrix
factorization of it, that means $\varphi\psi=\psi\varphi=f\cdot \mathbf 1_2$ and
$\mathrm{det}\,\varphi \cdot \mathrm{det}\,\psi=f^2$. Since $f$ is irreducible, we may consider
$\mathrm{det}\,\varphi=\mathrm{det}\,\psi=f$, and so, $\psi$ is the adjoint of $\varphi$. Therefore, it is sufficient to find
$\varphi= \left( \begin{smallmatrix}
   \varphi_{11} & \varphi_{12}
\\ \varphi_{21} & \varphi_{22}
         \end {smallmatrix}\right)$
such that $\mathrm{det}\,\varphi=f$ and $\varphi_{11}$ and $\varphi_{21}$ are two linearly independent linear forms.
 Applying some elementary transformations on the matrix  $\varphi$, we may suppose that:
\[
   \left \{\begin{array}{llll}
      \varphi_{11}=y_1-\lambda_1y_3 & \textrm{and} & \varphi_{21}=y_2-\lambda_2y_3, & \lambda_1,\lambda_2 \in k \\
                                    & \textrm{or} \\
      \varphi_{11}=y_1-\lambda y_2  & \textrm{and} & \varphi_{21}=y_3, & \lambda \in k.
   \end{array} \right.
\]
Let us consider the first case, when
$\varphi= \left( \begin{smallmatrix}
         y_1-\lambda_1y_3 & \varphi_{12}\\
         y_2-\lambda_2y_3 & \varphi_{22}
         \end{smallmatrix}\right)$.

Notice that $(\mathrm{det}\,\varphi)(\lambda_1,\lambda_2,1)=0$. Therefore $\lambda=(\lambda_1:\lambda_2:1)$ is a point on
the curve $V(f)$.
We want to show that $\varphi \sim \varphi_\lambda$.\\
For this, consider the product $\psi_\lambda \cdot \varphi$ that has the form
$\psi_\lambda \cdot\varphi= \left( \begin{smallmatrix}
                             f & g \\
                             0 & f
                           \end{smallmatrix}\right)$
with
$g=(\,y_1^2+(\lambda_1+1)y_1y_3+(\lambda_1^2+\lambda_1)y_3^2\,)\cdot\varphi_{12}-(\,y_2y_3+\lambda_2y_3^2\,)\cdot\varphi_{22}$.
Since $g\cdot(y_1-\lambda_1y_3)=\varphi_{12}\cdot
f-(y_2y_3+\lambda_2y_3^2)\cdot\mathrm{det}\,\varphi=f\cdot(\varphi_{12}-y_2y_3-\lambda_2y_3^2)$ and
$f$ is irreducible, we can write $g=f\cdot g_1$ with $g_1\in k[y_1,y_2,y_3]$. Therefore, we have
$\psi_\lambda\varphi= f\centerdot\left(\begin{smallmatrix}
                                 1 & g_1 \\
                                 0 & 1
                                 \end{smallmatrix}\right).$
Multiplying at left with $\varphi_\lambda$, we obtain
$\varphi= \varphi_\lambda\centerdot\left(\begin{smallmatrix}
                                    1 & g_1 \\
                                    0 & 1
                                  \end{smallmatrix}\right),$ that implies $\varphi \sim \varphi_\lambda$.

 The second case ($\varphi_{11}=y_1-\lambda y_2$ and $\varphi_{21}=y_3; \lambda \in k )$
can be treated exactly as above, replacing $\psi_\lambda$ with $\psi_{\lambda_0}$,
where $\lambda_0$ denotes the point (0:1:0).

(3) Because of the degrees of the entries, no module from  $\mathcal M_1 \cup \{\mathrm{Coker}\,\psi_s\}$ is isomorphic
with a module from $\mathcal M_{-1} \cup \{\mathrm{Coker}\,\varphi_s\}$.\\ Since any two equivalent matrices have
the same fitting ideals, for the rest, it is enough to consider the following fitting ideals:

$\bullet$ The modules from $\mathcal M_{-1} \cup \mathcal M_1$:\\
$\textrm{Fitt}_1(\varphi_\lambda)=\textrm{Fitt}_1(\psi_\lambda)=\langle
y_1-\lambda_1y_3,\,y_2-\lambda_2y_3,\,y_3^2 \rangle, \lambda=(\lambda_1:\lambda_2:1)\in V(f)$\\
 $\textrm{Fitt}_1(\varphi_{\lambda_0})=\textrm{Fitt}_1(\psi_{\lambda_0})=\langle y_1,\,y_3,\,y_2^2
\rangle, \lambda_0=(0:1:0)$.

$\bullet$ The modules from $\underline{\mathcal M}$:  $\textrm{Fitt}_1(\varphi_s)=\textrm{Fitt}_1(\psi_s)=\langle y_1,\,y_2 \rangle$.

(4) Follows from Corollary 6.4, \cite{Ei}.

(5) By construction, the modules of $\mathcal M_1$ are the syzygies of the modules of $\mathcal M_{-1}$.
Since
$\varphi_\lambda^t=\left(\begin{smallmatrix}
                  0 & 1 \\
                 -1 & 0
                 \end{smallmatrix}\right)
\psi_\lambda\left(\begin{smallmatrix}
              0 & -1 \\
              1 & 0
              \end{smallmatrix}\right)$, $\mathrm{Coker}\,\psi_\lambda\cong(\mathrm{Coker}\,\varphi_\lambda)^{\vee}$.
\end{pf}

\begin{thm}\label{deg1bundle}
Let $Y=\text{Proj}R$, that is the simple node singularity. Then:\\
(1) The coherent sheaves associated to the modules from $\mathcal M_{1}$ give all the isomorphism classes of line bundles of degree 1 over
$Y$;\\
(2) The coherent sheaves associated to the modules from $\mathcal M_{-1}$ give all the isomorphism classes of line bundles of degree -1 over
$Y$;\\
(3) The coherent sheaves associated to the modules from $\underline{\mathcal M}$ are not locally free.
\end{thm}

\begin{pf}
(3) The only singular point of  $V(f)$ is $(0:0:1)$. So, by \cite{Pfister}(1.3.8), it is sufficient to prove that $\textrm{Fitt}_1(\psi_s)R_{\langle y_1,y_2 \rangle}\ne
R_{\langle y_1,y_2 \rangle}$. Indeed, $\textrm{Fitt}_1(\psi_s)R_{\langle y_1,y_2 \rangle}=\langle y_1,y_2 \rangle R_{\langle y_1,y_2
\rangle}$, so $\Co\psi_s$ and its dual, $\Co\varphi_s$ are non-locally free.

(1) Any line bundle of degree one on $Y$ has the form $\mathcal O_Y(P)$, with $P$ regular point of $Y$.
 Following the proof of Theorem 3.8 from (\cite{LPP}), we obtain that the graded MCM $R$-module corresponding to
$\mathcal O_Y(P)$  is a module from $\mathcal M_1$, for any regular point $P$ of $Y$. 

(2) It follows from $\ref{2gen}$ and (1).
\end{pf}


 \subsection{Three-generated graded MCM R-modules }

For any $\lambda=(\lambda_1:\lambda_2:1)$ in $V(f)$ let be:
\[
\alpha_\lambda=\left(\begin{matrix}
                0 & y_1-\lambda_1y_3 & y_2-\lambda_2y_3
             \\y_1 & y_2+\lambda_2y_3 & (\lambda_1^2+\lambda_1)y_3
             \\y_3 &             0    & -y_1-(\lambda_1+1)y_3
                \end{matrix}\right)
\]
and $\beta_\lambda$ the adjoint of $\alpha_\lambda$.

Consider also the maps $\alpha_\lambda : R(-2)^3 \longrightarrow R(-1)^3$, given by the matrices $\alpha_\lambda$.\\

\begin{thm}  \label{3gen}
 For all $\lambda=(\lambda_1:\lambda_2:1)\in V(f), (\alpha_\lambda, \beta_\lambda)$ is a matrix factorization of $f$ and
 the set of three-generated graded MCM $R$-modules
\[\mathcal M_0=\{\,\mathrm{Coker}\,\alpha_\lambda\,|\,\lambda \in V(f)_{\mathrm{reg}}\backslash\{\lambda_0\}\,\}\]
has the
following properties: \\
(1) All the modules from $\mathcal M_0$ have rank 1.\\
(2) Every two different modules from $\mathcal M_0$ are not isomorphic.\\
(3) Every three--generated, rank 1, graded MCM $R$-module is isomorphic with one
of the modules from $\mathcal M_0$ or to $\mathrm{Coker}\,\alpha_s$.
\end{thm}

\begin{pf}
Clearly $\alpha_\lambda\beta_\lambda=\beta_\lambda\alpha_\lambda=f \cdot \mathbf 1_3$
for any $\lambda=(\lambda_1:\lambda_2:1) \in V(f)$.

(1) Since $\mathrm{det}\,(\alpha_\lambda)=f$, by Corollary 6.4 (\cite{Ei}), $\mathrm{Coker}\,\alpha_\lambda$ has rank 1.

(2) Suppose that there exist two invertible matrices, $U$ and $V$, with entries in $k$,  such that
$U\alpha_\lambda=\alpha_\xi V$ for $\lambda,\xi\in V(f)$. With the help of computer (we use \textsc{Singular}[GPS]) we obtain
that $\lambda=\xi$:
\begin{verbatim}
LIB "matrix.lib"; option(redSB);

ring r=0,(y(1..3),u(1..9),v(1..9),a,b,c,d),(c,dp);
ideal I=a3+a2-b2,c3+c2-d2;
qring Q=std(I);
matrix A[3][3]=   0, y(1)-a*y(3),      y(2)-b*y(3),
               y(1), y(2)+b*y(3),      (a2+a)*y(3),
               y(3),           0, -y(1)-(1+a)*y(3);
matrix B=subst(A,a,c,b,d);
matrix U[3][3]=u(1..9); matrix V[3][3]=v(1..9); int i;
matrix C=U*A-B*V;
ideal I=flatten(C);
ideal J=ideal(det(U)-1);
for (i=1;i<=3;i++)
 {J=J+transpose(coeffs(I,y(i)))[2];}
ideal L=std(J); L;
\end{verbatim}
 The first two entries of the ideal L are:
\begin{verbatim}
L[1]=b-d
L[2]=a-c
\end{verbatim}
Therefore $a=c$ and $b=d$, that means $\lambda=\xi$.

(3) Let be $M$ a three-generated, rank one, graded MCM $R$-module and
$(\varphi,\psi)$ the corresponding graded reduced matrix factorization. We can suppose $\mathrm{det}\,\varphi$=$f$ and
$\mathrm{det}\,\psi=f^2$. So, all entries of $\varphi$ have degree 1.
Since $f\in \langle y_1,y_3 \rangle$, by \cite{Eis2}, $\varphi$ has generalized zeros. Thus after some elementary
transformations, $\varphi=\left(\begin{smallmatrix}
          0& \varphi_1 & \varphi_2
         \\ \varphi_3&   a      & b
         \\ \varphi_4 &  c      & d
        \end{smallmatrix}\right)$
with  $\{\varphi_1, \varphi_2\}, \{\varphi_3, \varphi_4\}$ linearly
independent.\\
As $f\in \langle \varphi_1, \varphi_2 \rangle \cap \langle \varphi_3, \varphi_4 \rangle$, we can suppose
that $\varphi_1$ and $\varphi_3$ have non-zero coefficient of $y_1$. So, we can choose $\varphi_i, i=\overline{1,4}$ as
follows:
\[\left \{\begin{array}{lll}
   \varphi_1=y_1-\lambda_1y_3,\,\varphi_2=y_2-\lambda_2y_3 & \textrm{or} & \varphi_1= y_1-\lambda y_2,\,\varphi_2=y_3\\
   \varphi_{3}=y_1-\xi_1y_2,\,\varphi_{4}=y_2-\xi_2y_3 & \textrm{or} & \varphi_{3}=y_1-\xi y_2,\,\varphi_4=y_3.
   \end{array} \right.\]
Since $\mathrm{det}\,\varphi=f$, the points $(\lambda_1:\lambda_2:1),(\xi_1:\xi_2:1),(\lambda :1:0),(\xi :1:0)$ lay in
$V(f)$.
Therefore $\lambda=\xi=0$.\\
For any $\lambda=(\lambda_1:\lambda_2:1)$ in $V(f)$, we write  $\varphi_{1\lambda}=y_1-\lambda_1y_3$,
$\varphi_{2\lambda}=y_2-\lambda_2y_3$ and for $\lambda=(0:1:0)$ we write $\varphi_{1\lambda}=y_1,\varphi_{2\lambda}=y_3$.
Then $\varphi$ has the form:
 
 \hspace{2cm} $\varphi=\left(\begin{matrix}
        0 & \varphi_{1\lambda} & \varphi_{2\lambda}
       \\ \varphi_{1\xi} &   a & b
       \\ \varphi_{2\xi} &  c & d
       \end{matrix}\right)$ with $a,b,c,d$ linear forms.\\
Notice that, since $f \notin \langle y_1^2,y_1y_3,y_3^2\rangle$, it is not possible that $\lambda=\xi=(0:1:0)$.\\
To finish the proof, we need two helping results:
\begin{lem} \label{lema1}
    Let $M$ be a three-generated, rank one, graded MCM $R$-module and $(\varphi,\psi)$ a matrix factorization of $M$,
$\varphi$ having the above form. Then there exists $\lambda^\prime\in V(f)\backslash\{(0:1:0)\}, a^\prime, b^\prime,
c^\prime, d^\prime$
linear forms such that the matrix
 $$\varphi^\prime=\left(\begin{matrix}
                 0 & \varphi_{1\lambda ^\prime} & \varphi_{2\lambda^\prime}
                \\ y_1 & a^\prime & b^\prime
                \\ y_3 & c^\prime & d^\prime
                \end{matrix}\right)$$
together with its adjoint matrix $\psi^\prime$ form another matrix factorization $(\varphi^\prime,\psi^\prime)$ of $M$.
\end{lem}

\begin{pf}
We have to prove that after some elementary transformation the matrix $\varphi$ will become $\varphi^\prime$, that means, there exist two invertible $3 \times 3$ matrices $U$ and $V$, with entries in $k$, such that
$U\varphi^\prime \,=\,\varphi V$. For this, it is sufficient to prove that there exist two
invertible $3 \times 3$ matrices $U,V$ such that the first column of $U^{-1}\varphi V$ is
$\left( \begin{smallmatrix}
 0 \\
 y_1 \\
 y_3
\end{smallmatrix} \right)$.\\
Considering $U=(u_{ij})_{1\le i,j\le 3}$ and $V=(v_{ij})_{1\le i,j\le 3}$, the above condition lead to the following
system of equations:
\[
 \left\{\begin{array}{c}
 \varphi_{1\lambda}v_{21}+ \varphi_{2\lambda}v_{31}=y_1 u_{12}+ y_3u_{13}
 \\ \varphi_{1\xi}v_{11}+av_{21}+bv_{31}=y_1u_{22} +y_3u_{23}
 \\ \varphi_{2\xi}v_{11}+cv_{21}+dv_{31}=y_1u_{32}+y_3u_{33}.
 \end{array}\right.\]
In particular, $\varphi(0,1,0)\cdot\left(\begin{smallmatrix}
                                    v_{11} \\
                                    v_{21} \\
                                    v_{31}
                                    \end{smallmatrix}\right)=\left(\begin{smallmatrix}
                                                        0\\
                                                        0\\
                                                        0 \end{smallmatrix}\right)$.

Since $\mathrm{det}\,(\varphi(0,1,0))=f(0,1,0)=0$, we may choose a non-zero solution $(v_{11}, v_{21}, v_{31})$
which can be completed to an invertible matrix $V$ and
such that also the corresponding $(u_{12}, u_{13}, u_{22}, u_{23}, u_{32}, u_{33})$ can be completed to an invertible
matrix $U$.
\end{pf}

\begin{lem} \label{lema2}
  Let $M$ be a three-generated, rank one, graded MCM $R$-module and $(\varphi,\psi)$ a matrix factorization of $M$,
$\varphi$ having the form:
\[    \varphi=\left(\begin{array}{ccc}
              0 & \varphi_{1\lambda} & \varphi_{2\lambda}
              \\ y_1 &   a & b
               \\ y_3 &  c & d
             \end{array}\right)\] where $a,b,c,d$ are linear forms.
Then $(\alpha_\lambda,\beta_\lambda)$ is another matrix factorization of $M$.
\end{lem}
\begin{pf}\ \\
We make  some elementary transformations to simplify the entries $a,b$ and $c$. First, we eliminate the variable $y_1$ using the first column and the first line. Since $\lambda \ne (0:1:0)$, using the first line we can eliminate the variable $y_21$ from $b$. To "kill" the new $y_1$ in $a$, we subtract the first column from the second one. Therefore, instead of $b$ we can write $b y_3$ with $b\in k$.

\hspace{0,3cm}Consider the following polynomials (2-minors of $\alpha_\lambda$ and $\varphi$):

$\gamma=\left\vert\begin{matrix}
            y_1 & b y_3\\
            y_3 & d
          \end{matrix}\right\arrowvert$, 
$\delta=\left\vert\begin{matrix}
           y_1 & a \\
           y_3 & c
         \end{matrix}\right\vert$, 
$\bar{\gamma}=\left\vert\begin{matrix}
               y_1 & (\lambda_1^2+\lambda_1)y_3\\
               y_3 & -y_1-(\lambda_1+1)y_3
              \end{matrix} \right \vert$, 
$\bar{\delta}=\left\vert\begin{matrix}
               y_1 & y_2+\lambda_2y_3\\
               y_3 & 0
               \end{matrix}\right\vert$.

 Since
det$\varphi$ = det$\alpha_{\lambda}=f$, it holds the equality:
$\varphi_{1\lambda}(\bar{\gamma}-\gamma)=\varphi_{2\lambda}(\bar{\delta}-\delta).\hspace{0,2cm} (\ast)$\\
So $\varphi_{1\lambda}\mid \bar{\delta}-\delta$.
But $\bar{\delta}-\delta =-c(y_1-\lambda_1y_3)-y_3(y_2+\lambda_2y_3+\lambda_1 c-a)$ and
$a,c\in\langle y_2,y_3\rangle_k$. Therefore, $a=y_2+\lambda_2y_3+\lambda_1 c$ and $\bar{\delta}-\delta =-c
\varphi_{1\lambda}$.
Replacing  $\bar{\delta}-\delta$ in $(\ast)$ we get $\bar{\gamma}-\gamma=-c(y_2-\lambda_2y_3)$.
But $\bar{\gamma}-\gamma=y_1(-y_1-(\lambda_1+1)y_3)-d)-y_3^2(\lambda_1^2+\lambda_1-b)$ and
$c \in\langle y_2,y_3\rangle_k$. Therefore $d=-y_1-(\lambda_1+1)y_3, b=\lambda_1^2+\lambda_1$ and $c=0$. This shows that
$\varphi \sim \alpha_\lambda$.
\end{pf}
Using Lemma $\ref{lema1}$ and Lemma $\ref{lema2}$ the proof of the theorem $\ref{3gen}$ is finished.  
\end{pf}

\begin{thm}\label{deg0}
Let $Y$ be the projective cone over $R$.\\
(1) The coherent sheaves associated to the modules from $\mathcal M_0$ give all the isomorphism classes of 
line bundles of degree 0 over $Y$;\\
(2) The coherent sheaf associated to $\mathrm{Coker}\,\alpha_s$ is not locally free. 
Together with the sheafifications of the modules from $\underline{\mathcal M}$  give all the
isomorphism classes of rank one, non-locally free aCM sheaves.
\end{thm}
\begin{pf}
By computing $\textrm{Fitt}_2(\alpha_\lambda)R_{\langle y_1,y_2 \rangle}$ we find that
$\mathrm{Coker}\,\alpha_\lambda$ is locally free if and only if $\lambda$ is a regular point of $Y$.
In the previous subsection we have proved that the modules from $\km_{-1}$ and $\km_1$  have degree --1, respectively 1.
After some shiftings, they give all line  bundles of degree $3t-1$ and $3t+1$, with $t\in \Z$.
The remaining rank one graded MCM $R$--modules (from $\km_0$), define  the line bundles of degree $3t$.\\
 1) We prove first that the sheafification of $\mathrm{Coker}\,\alpha_\xi$ has degree 0, where $\xi=(-1:0:1)$.
For this, we compute $(\textrm{Coker}\,\alpha_\xi \otimes \textrm{Coker}\,\alpha_\xi)^{\vee\vee}$, using the computer. The following two procedures were used in \cite{LPP}.
\begin{verbatim}
LIB "matrix.lib"; option(redSB);

proc reflexivHull(matrix M)
{ module N=mres(transpose(M),3)[3];
  N=prune(transpose(N));
  return(matrix(N));}

proc tensorCM(matrix Phi, matrix Psi)
{ int s=nrows(Phi); int q=nrows(Psi);
  matrix A=tensor(unitmat(s),Psi); 
  matrix B=tensor(Phi,unitmat(q));
  matrix R=concat(A,B);
  return(reflexivHull(R));}

ring R1=0,(y(1..3)),(c,dp);
ideal i=y(1)^3+y(1)^2*y(3)-y(2)^2*y(3);
qring S1=std(i);
matrix M[3][3]=   0, y(1)+y(3),  y(2),
               y(1),      y(2),     0,
               y(3),         0, -y(1);
tensorCM(M,M);
_[1,1]=0
\end{verbatim}
This means that $\textrm{Coker}\,\alpha_\xi$ is a self dual module, so the matrix $\alpha_\xi$ corresponds to $\mathcal B(0,1,-1)$.
The graded map $\alpha_\xi : R(-2)^3 \longrightarrow R(-1)^3$, should have therefore the degree of the form 3$t$.
From the graded exact sequence
\[0 \longrightarrow (\textrm{Coker}\,\alpha_\xi)^{\vee} \longrightarrow R(1)^3 \longrightarrow R(2)^3 \longrightarrow \textrm{Coker}\,\alpha^t_\xi \longrightarrow 0\]
we see that the degree of $\textrm{Coker}\,\alpha^t_\xi$ is $9-3t$. But there exists a graded isomorphism between $\textrm{Coker}\,\alpha_\xi$ and
 $\textrm{Coker}\,\alpha^t_\xi \otimes R(-3)$. So, $(9-3t)-9=3t$, that means $t=0$.\\
2) Let us now consider $\mathrm{Coker}\,\alpha_\lambda$, with $\lambda=(a:b:1)$ an arbitrary regular point on the nodal curve.
Its sheafification has degree of the form $3t$.

 Consider also the module $\textrm{Coker}\,\psi_\xi$, with the corresponding graded map\\
 $\psi_\xi: R(-2)^2 \To R \oplus R(-1)$. As we have seen in the previous subsection, this module has the degree 1.
The bundle corresponding to $\textrm{Coker}\,\psi_\xi$ has the form $\kb(1,1,\mu)$, the one corresponding to
$\mathrm{Coker}\,\alpha_\lambda$ has the form $\kb(3t,1,\mu')$.

By Serre duality, for any two vector bundles $\kf, \ke$ on the simple node, \[\textrm{Hom}(\kf,\ke)=\textrm{Ext}(\ke,\kf)=
\textrm{H}^1(\ke^{\vee}\otimes\kf).\]
So, $\textrm{Ext}^1(\kb(1,1,\mu),\kb(3t,1,\mu'))=\textrm{H}^1(\kb(-1,1,\mu^{-1})\otimes \kb(3t,1,\mu'))=\\
 \textrm{H}^1(\kb(3t-1,1,\mu^{-1}\mu'))$.
 But $\textrm{dim}_k(\textrm{H}^1(\kb(3t-1,1,\mu^{-1}\mu')))=1$ iff $3t-1=-1$, that means $t=0$ ( \cite{igor}, Kapitel 3, 4.1).

We prove that $\textrm{dim}_k(\textrm{Ext}^1(\textrm{Coker}\,\psi_\xi,\mathrm{Coker}\,\alpha_\lambda)=1$, for any
 $\lambda=(a:b:1)$ regular point of the simple node. With this we are done.\\
According to the theorem $\ref{mfext}$, a module $M$ with a graded extension
$$0\To \mathrm{Coker}\,\alpha_\lambda \To M \To \textrm{Coker}\,\psi_\xi \To 0,$$
is the cokernel of a graded map $T : R^5(-2)\To R^3(-1)\oplus R \oplus R(-1)$, given by a
square $5 \times 5$ matrix of the form
$\left(\begin{smallmatrix}
\alpha_\lambda & D\\
\large{0}&\psi_\xi
\end{smallmatrix}\right).$
$D$ is a $3 \times 2$ matrix with linear entries, of the form $D=\left(\begin{smallmatrix}
d_1& d_2\\
d_3&d_4\\
d_5&d_6
\end{smallmatrix}\right)$, such that $\beta_\lambda D \varphi_\xi=0$ in $R$.\\
We make some linear transformations on $T$, in order to obtain a simple form of $D$.
First of all, we can eliminate the variable $y_1$ in all entries of $D$, by subtracting one of the first three columns
multiplied with some constant. By subtracting the last line from the first one we eliminate the variable
$y_2$ in the entry $d_1$. We add in this way $y_1$ to the entry $d_2$, but we can "kill" it using the second column.
 In the same way, using the last line and the third column of $T$, we eliminate $y_2$ also in the entry $d_5$.
By subtracting the first column from the last one, we eliminate  the variable $y_3$ in $d_6$. We "kill" the new $y_1$
in the entry $d_4$ using the last line.\\
Let us  now study the relation $\beta_\lambda
D \varphi_\xi=0$ in $R$. For simplicity, we use for this the computer.

The procedure \texttt {condext} returns the ideal given by the coefficients of $y_1$ in the entries of a matrix, after
it reduces the entries to the polynomial $y^3_1+y_1^2y_3-y^2_2y_3$. This procedure will be used also in the last section.
\begin{verbatim}
 LIB"matrix.lib"; LIB"homolog.lib"; LIB"linalg.lib";

 proc simple(ideal P)
 { int j,i; poly F;
 list L=0;
 for(j=1;j<=size(P);j++)
    { L=factorize(P[j]);
      if(size(L[1])>2)
          {F=1;
           for (i=2;i<=size(L[1]);i++)
               {if (L[1][i]==y(2) or L[1][i]==y(3))
                     { L[1][i]=1;}
                F=F*L[1][i]^(L[2][i]);}  
           P[j]=F;}}
  return(P);}

 proc condext(matrix A,B,D)
 { matrix Aa=adjoint(A); matrix Ba=adjoint(B); matrix G=Aa*D*Ba;
  ideal g=flatten(G);
  matrix V; int k,j; ideal P=0; list L=0;
  for(j=1;j<=size(G);j++)
   { g[j]=reduce(g[j],std(y(1)^3+y(1)^2*y(3)-y(2)^2*y(3)));
     V=coef(g[j],y(1));
     for(k=1;k<=1/2*size(V);k++)
     { P=P+V[2,k];}}
  P=interred(P); P=simple(P); 
  return(P);}
\end{verbatim}
We define the ring $R$ and the matrices $\psi_\xi$, $\phi_\xi$, respectively $\alpha_\lambda$.
\begin{verbatim}
  ring R=0,(y(1..3),d(1..6),a,b),(c,dp(3),dp(6),dp(2));
  ideal i=y(1)^3+y(1)^2*y(3)-y(2)^2*y(3),a3+a2-b2;
  qring S=std(i); ideal P;
  matrix psi[2][2]= y(1)^2,-y(2)*y(3),
                     -y(2), y(1)+y(3);
  matrix phi[2][2]=y(1)+y(3), y(2)*y(3),
                        y(2), y(1)^2;
  matrix A[3][3]=0,y(1)-a*y(3),     y(2)-b*y(3),
              y(1),y(2)+b*y(3),     (a2+a)*y(3),
              y(3),          0,-y(1)-(a+1)*y(3);
\end{verbatim}
We define the matrix $D$ and put the condition $\beta_\lambda D \varphi_\xi=0$.
\begin{verbatim}
  matrix D[3][2]=d(1..6);
  P=condext(A,psi,D); P;
P[1]=y(2)*d(6)-y(3)*d(3)-y(3)*d(5)*a
P[2]=y(2)*d(2)+y(2)*d(5)+y(3)*d(1)*a+y(3)*d(1)-y(3)*d(5)*b
P[3]=y(2)*d(1)-y(2)*d(4)+y(3)*d(1)*b-y(3)*d(3)-y(3)*d(5)*a^2-
     -y(3)*d(5)*a
P[4]=d(1)*b-d(2)*a^2-d(4)*b-d(5)*a^2-d(6)*b
P[5]=-d(1)*a-d(1)+d(2)*b+d(4)*a+d(4)+d(5)*b+d(6)*a+d(6)
\end{verbatim}
From P[1], we obtain $y_3 | d_6$, but we have already eliminated the variable $y_3$ in $d_6$, so we can suppose
that $d_6=0$. Under this condition, P[1] implies that $d_3=-ad_5$.\\ From P[2] we obtain that $y_2 | d_1(a+1)-d_5b$, but $d_1$ and $d_5$ have no $y_2$.
Therefore, $d_1(a+1)-d_5b=d_2+d_5=0$. Since $a \neq 0$,  P[5] implies that $d_4=d_1$ and P[4] implies $d_5=d_1b/a^2$.
Since  $d_1$ has the form $d_1=a_1y_3$, with $a_1 \in k$ and we want a nonzero extension, we can choose $a_1=a^2$. \\The matrix $T$ becomes:
$$ T=\left(\begin{matrix}
0&y_1-ay_3&y_2-by_3&a^2y_3&-by_3\\
y_1&y_2+by_3&(a^2+a)y_3&-by_3&a^2y_3\\
y_3&0&-y_1-(a+1)y_3&by_3&0\\
0&0&0&y_1^2&-y_2y_3\\
0&0&0&-y_2&y_1+y_3
\end{matrix}\right).$$
Since there are no matrices $U,V$ such that $D=\alpha_\lambda U+V\psi_\xi$ the extension defined by $T$ is nonzero.\\
This proves that $\textrm{dim}_k(\textrm{Ext}^1(\textrm{Coker}\,\psi_\xi,\mathrm{Coker}\,\alpha_\lambda)=1$,
and with this we are done.
\end{pf}

\begin{cor}\label{3gen2}
Every three-generated, rank 2, indecomposable, graded MCM $R$-module is isomorphic to one
of the modules $\mathrm{Coker}\,\beta_\lambda$, $\lambda=(\lambda_1:\lambda_2:1)$ a point on $V(f)$.
\end{cor}

\begin{pf}
 Let $M$ be a three-generated, rank two, indecomposable, graded MCM $R$-module and $(\varphi,\psi)$ a matrix
factorization of $M$. Then $\mathrm{Coker}\,\psi$ is a three-generated, rank one, graded MCM $R$-module.
Therefore, it is isomorphic to one of the modules from $\mathcal M_0 \cup \{\mathrm{Coker}\,\alpha_s\}$ and
$\mathrm{Coker}\,\varphi$ is isomorphic to one of the modules from $\{\mathrm{Coker}\,\beta_\lambda|
\lambda=(\lambda_1:\lambda_2:1), \lambda\in V(f)\}$.
\end{pf}

We present a short overview over the results concerning the rank one graded MCM $R$--modules:

 The rank one graded MCM R-modules corresponding to the \textbf{line bundles} on $Y$ are given by the cokernel of some graded maps
defined by the following matrices, with $(\lambda_1:\lambda_2:1)$ a regular point on $Y$:
$$\left(\begin{matrix}
                           y_1-\lambda_1y_3 & y_2y_3+\lambda_2y_3^2
                        \\y_2-\lambda_2y_3 &y_1^2+(\lambda_1+1)y_1y_3+(\lambda_1^2+\lambda_1)y_3^2
                              \end{matrix}\right),
    \left(\begin{matrix}
           y_1+y_3 & y_2^2
            \\ y_3 & y_1^2
     \end {matrix}\right),$$
$$\left(\begin{matrix}
   y_1^2+(\lambda_1+1)y_1y_3+(\lambda_1^2+\lambda_1)y_3^2 & -y_2y_3-\lambda_2y_3^2
\\      -(y_2-\lambda_2y_3)                               &     y_1-\lambda_1y_3
                    \end{matrix}\right),\left(\begin{matrix}
                                        y_1^2 & -y_2^2
                                         \\ -y_3 & y_1+y_3
                                         \end {matrix}\right),$$
 $$\left(\begin{matrix}
    0&y_1-\lambda_1y_3 & y_2-\lambda_2y_3
\\y_1&y_2+\lambda_2y_3 & (\lambda_1^2+\lambda_1)y_3
\\y_3&      0    & -y_1-(\lambda_1+1)y_3
\end{matrix}\right).$$
The rank one graded MCM $R$-modules corresponding to the \textbf{non-locally free} aCM sheaves are given by the following set
of matrices:
$$\{\hspace{0,1cm} \left(\begin{matrix}
                                  y_1 & y_2y_3
                              \\  y_2 & y_1^2+y_1y_3
                                 \end{matrix}\right) ,
                                 \left(\begin{matrix}
                                  y_1^2+y_1y_3 &-y_2y_3
                                  \\-y_2 & y_1
                                  \end{matrix}\right),
                                  \left(\begin{matrix}
                                  0 & y_1 & y_2
                                  \\y_1 & y_2 & 0
                               \\y_3 & 0   & -y_1-y_3
                             \end{matrix}\right)    \hspace{0,1cm}\}.$$


\section { Rank two, graded,  MCM modules }

In the following there are described all isomorphism classes of rank two, in-\\decomposable MCM modules over
$R=k[y_1,y_2,y_3]/(f), f=y_1^3 +y_1^2y_3-y_2^2y_3$.\ \ \ \ \\
 In the case of locally free MCM modules, we use the classification of vector bundles over the projective cone
$Y=\textrm{Proj}\ R$ (see \cite{DG}).\\
 The matrix factorizations of the non--locally free, rank two, MCM modules are computed using Proposition $\ref{mfext}$.

\subsection{The classification of rank two locally free MCM R--modules }

There are two types of rank two, indecomposable, vector bundles on Proj$R$:\\
 $\bullet$ $\mathcal B(a,2,\lambda)$, with $a\in \mathbb Z $ and $\lambda \in k^*$ \\
 $\bullet$$\mathcal B(\mathbf{d},1,\lambda)$, with $\mathbf{d}$ a 2-cycle with entries in $\mathbb Z$ and $\lambda \in
k^*$($\mathbf{d}=(a,b), a\neq b$).

To generate the first type of rank two vector bundles it is sufficient to know the bundle $\mathcal B(0,2,1)$ and
the line bundles, because, for any $\lambda \in k^*$, \\ $$\mathcal B(a,2,\lambda)\cong \mathcal B(a,1,\lambda)\otimes \mathcal B(0,2,1).$$
The fact that the bundle $\mathcal B(0,2,1)$ is uniquely determined by the exact sequence
\begin{equation}
0 \longrightarrow \mathcal{O}_Y \longrightarrow \mathcal B(0,2,1) \longrightarrow \mathcal{O}_Y\longrightarrow 0,
  \end{equation}
  provide a way to determine the graded MCM $R$-module corresponding to it. Using the matrix factorizations of the rank 1 MCM modules 
 and a \textsc{Singular}--procedure that computes the tensor product of two locally free MCM modules, one can  
construct matrix factorizations for all bundles $\mathcal B(a,2,\lambda)$.

 The second type of rank two vector bundles can be generated by the bundles $\mathcal B((0,n),1,\lambda)$ and the line
bundles, using the tensor product formula:\\
  $\mathcal B((a,b),1,\lambda)\otimes \mathcal B(c,1,\mu)\cong \mathcal B((a+c,b+c),1,\lambda\mu^2)$.

  For $\mathbf{d}=(a,b)$ and $\mathbf{e}=(c,d)$ two 2-cycles with entries in $\mathbb Z$, we have:
 \[\mathcal B(\mathbf{d},1,\lambda)\otimes \mathcal B(\mathbf{e},1,\mu)\cong
  \mathcal B(\mathbf{f}_1,1,\lambda\cdot\mu)\oplus \mathcal B(\mathbf{f}_2,1,\lambda\cdot\mu),  \]
  where
  $\mathbf{f}_1=(a+c,b+d)$ and $\mathbf{f}_2=(a+d,b+c)$. If $\mathbf{f}_i=(\alpha,\alpha)$($i=1$ or 2), then
  $\mathcal B(\mathbf{f}_i,1,\lambda\cdot\mu)$ splits as: $\mathcal B(\mathbf{f}_i,1,\lambda\cdot\mu)=
   \mathcal B(\alpha,1,\sqrt{\lambda\cdot\mu})\oplus \mathcal B(\alpha,1,-\sqrt{\lambda\cdot\mu})$.

Therefore, inductively, we can obtain all $\mathcal B((0,n),1,\lambda)$, $n\in\mathbb{N}^*$, if we know the
bundles $\mathcal B((0,1),1,\lambda)$. By duality, ($\mathcal B(\mathbf{d},1,\lambda)^\vee\cong
\mathcal B(-\mathbf{d},1,\lambda^{-1})$) we obtain also $\mathcal B((0,n),1,\lambda)$ with n negative integer.\\
The bundles $\mathcal B((0,1),1,\lambda)$ are uniquely determined by the existence of the
exact sequences
\begin{equation}
0 \longrightarrow \mathcal{O}_Y \longrightarrow \mathcal B((0,1),1,\lambda) \longrightarrow \mathcal B(1,1,-\lambda)
\longrightarrow 0.   \end{equation}
 Using this we can compute the graded MCM $R$-modules corresponding to them.
So, inductively, one can obtain all rank two graded indecomposable locally free MCM $R$-modules.

 In the sequel we determine the module $M_2$ corresponding to $\mathcal B(0,2,1)$.

\begin{lem} \label{psi}
Let be $\rho=\left(\begin{matrix}
  y_1^2+y_1y_3 & -y_2&-y_3&0
\\-y_2y_3 & y_1& 0 & -y_3\\
    0    &  0 &y_1 & y_2
\end{matrix}\right)$,

$\psi=\left(\begin{matrix}
  y_1 & y_2&y_3&\phantom{-}0
\\y_2y_3 & y_1^2+y_1y_3& 0 & y_3\\
    0    &  0 & y_1^2+y_1y_3 & -y_2\\
    0    &  0 & -y_2y_3& y_1
\end{matrix}\right)$,
$\gamma=\left(\begin{matrix}
    0    &   0 & y_2y_3&y_1^2+y_1y_3
\end{matrix}\right)$ \\and
$\varphi= \left(\begin{matrix}
\rho\\ \gamma
\end{matrix}\right)$.
Then $(\psi,\varphi)$ is a matrix factorization of $\Omega^1 _R(m)$, where $m$ is the unique graded
maximal ideal of $R$, $m=\langle y_1,y_2,y_3\rangle$. More, the following exact sequence
\begin{equation}
 \begin{CD}
 \stackrel{\psi}{\longrightarrow} R(-3)\oplus R(-2)^3 \stackrel{\rho}{\longrightarrow}
R(-1)^3@>(y_1,y_2,y_3)>>\textstyle{m}
 \longrightarrow 0
 \end{CD}
 \end{equation}
  is a graded minimal free resolution of $m$. In particular, $\Omega^1 _R(m)$ has no free summands. \end{lem}
 \begin{pf}
Clearly $\varphi\psi=\psi\varphi=f\cdot \mathbf 1_4$ and the above sequence is a complex.
Let $u_1,u_2,u_3 \in R$ such that $\sum_{i=1}^3 y_i u_i=0$. We show that
$u=\left(\begin{smallmatrix}
    u_1 \\u_2 \\u_3
\end{smallmatrix}\right)$ is an element of $\textrm{Im}\rho$.
Subtracting multiples of the second and third columns of $\rho$ from $u$ we may suppose that $u_1$ depends only on $y_1$.
As the maps are graded, we may suppose that
$u$ is graded, so $u_1=ay_1^s, a\in k, s\in \mathbb N$.\\
If $a\neq 0$, $s+1 \geq 3$, so, subtracting from $u$ multiples of
$\left(\begin{smallmatrix}
    y_1^2 \\-y_2y_3\\y_1^2
\end{smallmatrix}\right)$ $\in \textrm{Im}\rho$, we reduce to the case $u_1=0$. Then $y_2u_2+y_3u_3=0$ and, since
$\{y_2,y_3\}$ is a regular sequence in $R$, $u$ is a multiple of the fourth column of $\rho$.\\ To show that
$\textrm{Ker}\,\rho \subset \textrm{Im}\,\psi$, it is enough to prove that $\textrm{Ker}\,\rho \subset
\textrm{Ker}\,\varphi$.\\
Denote by $\rho_3$ the third row of $\rho$ and choose $\nu$ an element of $\textrm{Ker}\,\rho$.
Since $y_1\gamma=y_2y_3\rho_3$ and
 $\rho_3\nu=0$ we get $\gamma\nu=0$ and $\nu $ is an element of $\textrm{Ker}\,\varphi$.\\
 Because no entry of $\varphi$ or $\psi$ is unite, $\Omega^1 _R(m)$ has no free summands.
\end{pf}

 \begin{prop}
 There exists a graded exact sequence:
 \begin{equation}
 0 \longrightarrow R \stackrel{i}{\longrightarrow} \Omega^2 _R(m)\otimes R(3)\stackrel{\pi}{\longrightarrow} m
\longrightarrow 0 \end{equation}
and
$\Omega^2 _R (m)\otimes R(3)$ corresponds to the bundle $\mathcal B(0,2,1)$.
\end{prop}

 \begin{pf}
 1) We prove the existence of the exact sequence (4).\\
 Define the map $i:R\longrightarrow \Omega^2 _R (m)\otimes R(3)$ by $i(1)=\left(\begin{smallmatrix}
   0 \\y_3\\-y_2\\y_1
\end{smallmatrix}\right)$ (it is the fourth column of $\psi$) and let $\pi: \Omega^2 _R(m)\otimes R(3)\longrightarrow m$ be the
projection on the first component. Since $\Omega^2 _R(m)\otimes
R(3)=\textrm{Im}\,\psi\otimes R(3)\subset R\oplus R(1)^3$, $i$ and $\pi$ are graded morphisms. Clearly, $i$ is injective,
$\pi$ is surjective and
$\mathrm{Im}\, i\subset \mathrm{Ker}\, \pi$. We prove that
$\mathrm{Ker}\,\pi\subset \textrm{Im}\,i$.\\
Let $v=\psi\cdot \left(\begin{smallmatrix}
    a \\b\\c\\d
\end{smallmatrix}\right)$ be an element in $\textrm{Ker}\,\pi$. Then $\left(\begin{smallmatrix}
    a \\b\\c
\end{smallmatrix}\right)$$\in \textrm{Im}\,\rho$. \\
Denote by $\psi '$ the $4\times 3 $ matrix obtained from $\psi$ by
eliminating the last column. \\
Then $v=\psi'\cdot \left(\begin{smallmatrix}
    a \\b\\c
\end{smallmatrix}\right)+d\left(\begin{smallmatrix}
   0 \\y_3\\-y_2\\y_1
\end{smallmatrix}\right).$ Since the columns of $\psi'\rho$ are in $ R\cdot \left(\begin{smallmatrix}
    0 \\y_3\\-y_2\\y_1
\end{smallmatrix}\right)$, $v$ is in $\textrm{Im}\,i$.

2) We prove the  indecomposability of the $R$-module $M_2=\Omega^2 _R(m)\otimes R(3)$.\\
If it would decompose, it would be isomorphic to $\textrm{Coker}\,\theta$, with $\theta$ of the
form $\left(\begin{smallmatrix}     A&0\\
    0&B
\end{smallmatrix}\right)$, $A$ and $B$ quadratic matrices with determinant equal to $f$.
They define rank 1, graded MCM $R$-modules, so, they are equivalent to one of the matrices
$\{\varphi_\lambda, \psi_\lambda |\, \lambda\in V(f)\}$. Since $\theta\sim\varphi$, $\textrm{Fitt}_2(\theta)= \textrm{Fitt}_2(\varphi)=m^2$.\\
But the elements of degree 2 from
$\textrm{Fitt}_2(\theta)$ are given just by $l^A _1\cdot l^B _1$, $l^A _1\cdot l^B _2,$ $l^A _2\cdot l^B _1,$ $l^A
_2\cdot l^B _2$ where  $l^A _1, l^A _2$ respectively $ l^B _1, l^B _2$ are the entries of  $A$ and $B$ of degree 1.
The ideal generated by them can not be $m^2$, that is minimally generated by 6 elements. Therefore, $M_2$ is
indecomposable.

3) To prove that $M_2$ is locally free, it is sufficient to notice that\\
 $\textrm{Fitt}_2 (\varphi)R_{\langle y_1,y_2\rangle}=R_{\langle y_1,y_2\rangle}$ and $\textrm{Fitt}_3
(\varphi)=0$.

4) From the exact sequence $\textstyle(4)$ we get the following  exact sequence of vector bundles on $Y=\mbox{Proj}R$:
\[0 \longrightarrow \mathcal O_Y \longrightarrow \widetilde{M_2} \longrightarrow \mathcal O_Y
\longrightarrow 0.  \]
$\widetilde{M_2}$ is an indecomposable vector bundle of rank 2, so it is isomorphic to $\mathcal B(0,2,1)$.
\end{pf}

\begin{thm}
 The rank two vector bundles of type $\mathcal B(a,2,\lambda)$, with $\lambda\in k^*$ and $a$ integer, are sheafification of
 the $R$-modules ${(M_2 \otimes L)}^{\vee\vee}\otimes R(k)$, with $L$ in $\mathcal{M}_0\cup\km_1\cup \km_{-1}$ and $k \in \Z$.
 More, the modules corresponding to
   \[\left\{\begin{array}{ll}
 \mathcal B(0,2,\lambda), \lambda\neq 1, & \phantom{-}\, \textrm{have 6 generators};\\
 \mathcal B(-1,2,\lambda), & \phantom{-}\,  \textrm{have 4 generators};\\
 \mathcal B(1,2,\lambda), & \phantom{-}\,   \textrm{have 4 generators}.\end{array}\right.
\]
\end{thm}

\begin{pf}
The first statement is follows directly from the previous proposition.\\
 The modules corresponding to the bundles $\mathcal B(0,2,\lambda)$ are
$(M_2\otimes\textrm{Coker}\,\alpha_\lambda)^{\vee\vee}$ and
we can compute their matrix factorizations using the computer. 
\begin{verbatim}
proc M2(ideal I)
{ matrix A=syz(transpose(mres(I,3)[3]));
  return(transpose(A));}

setring S1;
ideal I=maxideal(1); matrix C=M2(I);
ring R2=0,(y(1..3),a,b),(c,dp);
ideal I=y(1)^3+y(1)^2*y(3)-y(2)^2*y(3),a3+a2-b2;
qring S2=std(I);
matrix A[3][3]=    0, y(1)-a*y(3),      y(2)-b*y(3),
                y(1), y(2)+b*y(3),      (a2+a)*y(3),
                y(3),           0, -y(1)-(a+1)*y(3);
matrix C=imap(S1,C);
\end{verbatim}
We compute the dimension of the matrix corresponding to $(M_2\otimes \textrm{Coker}\,\alpha_\lambda)^{\vee\vee}$.
\begin{verbatim}
nrows(tensorCM(C,A));\\tensorCM is defined in the previous section
6
\end{verbatim}
The same, the modules corresponding to  $\mathcal B(-1,2,\lambda)$ are
 $(M_2\otimes\textrm{Coker}\,\varphi_\lambda)^{\vee\vee}$. We compute the size of the matrix corresponding to this
module.
\begin{verbatim}
matrix B[2][2]=y(1)-a*y(3),                  y(2)*y(3)+b*y(3)^2,
               y(2)-b*y(3),y(1)^2+(a+1)*y(1)*y(3)+(a2+a)*y(3)^2;
nrows(tensorCM(C,B));
4
B= y(1)+y(3),y(2)^2,
   y(3),       y(1)^2;
nrows(tensorCM(C,B));
4
\end{verbatim}
 By duality, it follows also the last statement.  
 \end{pf}

Let $\xi=(-1:0:1)$ and $\lambda_0=(0:1:0)$ two regular points on  $Y$. We consider the graded
maps $\psi_\xi : R(-2)^2 \longrightarrow R \oplus R(-1)$ and  $\alpha_{\lambda_0} : R(-2)^3 \longrightarrow R(-1)^3$,
given
by the matrices with the same name. (see the previous section)\\

 \begin{lem} \label{degree}
If we denote by $\mu_0$ the regular point on the nodal curve $Y$ such that
$\widetilde{\mathrm{Coker}\,\psi_\xi}\cong \mathcal B(1,1,\mu_0)$, then $\mathcal O_Y(1)=\mathcal B(3,1,-\mu_0^3)$
and the line bundle $\mathcal B(1,1,-\mu_0)$ is given by $\mathrm{Coker}\,\psi_{\lambda_0}$.
\end{lem}

\begin{pf}
 By $\ref{deg1bundle}$, the degree of $\textrm{Coker}\,\psi_{\lambda_0}$ is 1. 
We denote $\widetilde{\textrm{Coker}\,\psi_{\lambda_0}}$ by $\mathcal B(1,1,\theta)$.\\
Let us compute, the tensor product ($\textrm{Coker}\,\alpha_\xi \otimes \textrm{Coker}\,\psi_\xi)^{\vee\vee}$.
 \begin{verbatim}
matrix N[2][2]= y(1)^2, -y(2)*y(3),
                  -y(2),  y(1)-y(3);
matrix L=tensorCM(M,N);
print(L);
  y(1),            -y(3),
y(2)^2, -y(1)^2+y(1)*y(3)
\end{verbatim}
 After some elementary transformations, the matrix L becomes $\psi_{\lambda_0}$, so we have
obtained that ($\textrm{Coker}\,\alpha_\xi \otimes \textrm{Coker}\,\psi_\xi)^{\vee\vee}\cong
\textrm{Coker}\,\psi_{\lambda_0}$.
By sheafification, it becomes
$\mathcal B(0,1,-1)\otimes\mathcal B(1,1,\theta)\cong\mathcal B(1,1,-\mu_0)$, therefore $\theta=\mu_0$.

Let us now compute the reflexive hull of
$\textrm{Coker}\,\psi_{\lambda_0}\otimes\textrm{Coker}\,\psi_{\lambda_0}\otimes\textrm{Coker}\,\psi_{\lambda_0}$.
\begin{verbatim}
matrix L[2][2]= y(1)^2, -(y(1)^2+y(2)^2),
                 -y(3),             y(1);
tensorCM(L,tensorCM(L,L));
_[1,1]=0
\end{verbatim}
This means that
($\textrm{Coker}\,\psi_{\lambda_0}\otimes\textrm{Coker}\,\psi_{\lambda_0}\otimes\textrm{Coker}\,\psi_{\lambda_0})^{\vee}\cong R$,
so the line bundle $\mathcal O_Y(1)$ is $\mathcal B(3,1,-\mu_0^3).$    \end{pf}

Let us now determine the graded MCM modules corresponding to the bundles $\mathcal B((0,1),1,\lambda)$ with
$\lambda \in k^*$.
Consider the module $\Omega^1_R(M_2)=\mathrm{Coker}\,\psi$, where \\$\psi : R(-4)^3 \oplus R(-3) \longrightarrow R(-3)
\oplus R(-2)^3$ is given as in lemma $\ref{psi}$.\\

\begin{lem} \label{generatori}
Consider $\lambda$ a regular point on the curve.
\begin{enumerate}
\item $(\Omega^1_R(M_2)\otimes\textrm{Coker}\,\alpha_\lambda)^{\vee\vee}$ has
 \[
   \left\{\begin{array}{ll}
   \textrm{4 generators,}&\textrm{if $\lambda=(1:0:1)$};\\
   \textrm{3 generators,}&\textrm {otherwise;}
          \end{array}\right.
 \]
\item $(\Omega^1_R(M_2)\otimes\textrm{Coker}\,\varphi_\lambda)^{\vee\vee}$ has 5 generators;
\item $(\Omega^1_R(M_2)\otimes\textrm{Coker}\,\psi_\lambda)^{\vee\vee}$ has 5 generators;
\end{enumerate}
\end{lem} 
 \begin{pf}
 We define the module $\Omega^1_R(M_2)$ by:
\begin{verbatim}
matrix D=transpose(syz(C));
\end{verbatim}
and use the procedure \texttt{tensorCM} as before.   \end{pf}

 We compute the matrix corresponding to the MCM module\\ $(\Omega^1 _R(M_2)\otimes
\textrm{Coker}\varphi_\xi)^{\vee \vee}$.
\begin{verbatim}
matrix D=transpose(syz(C));
matrix B[2][2]=y(1)+y(3), y(2)*y(3),
                    y(2),    y(1)^2;
matrix A=tensorCM(D,B);
\end{verbatim}
After some linear transformations, the matrix $A$ it becomes:
$$A= \left(
\setlength{\unitlength}{0,3cm}
\begin{picture}(10,6)
\thinlines  \put(-0.2,0){\line(9,0){10}}
\thinlines  \put(4,6){\line(0,-2){11}}
\put (6,5){$y_1$}
\put (8.5,5){0}
\put (6,3){0}
\put (8.5,3){0}
\put (6,1){0}
\put (8,1){$-y_3$}
\put (1,3){\large $\alpha_{\xi}$}
\put (1,-4){\huge 0}
\put (6,-3){\large $\psi_\xi$}
\end{picture}
\right).$$
Observation: This matrix is linear equivalent to the matrix $T$ obtained in the proof of Theorem $\ref{deg0}$,
for $a=-1$ and $b=0$.\\

\begin{prop} \label{omega}
\begin{enumerate}
\item The graded module corresponding to $\mathcal B((0,1),1,\mu_0)$ is $(\Omega^1 _R
(M_2)\otimes\textrm{Coker}\,\varphi_\xi)^{\vee\vee}\otimes R(2)$.
\item $\widetilde{\Omega^1_R(M_2)}=\mathcal B((-4,-5),1,\mu_0^{-9})$.
\end{enumerate}
\end{prop}

\begin{pf}
(1) Since $\beta_{\xi}\left( \begin{smallmatrix}y_1&0\\0&0\\0&-y_3\end{smallmatrix}\right)\varphi_{\xi}=0$ in $R$, by $\ref{mfext1}$, there exists the graded extension 
\[
0 \longrightarrow \textrm{Coker}\,\alpha_\xi \longrightarrow \textrm{Coker}\,A \longrightarrow \textrm{Coker}\,\psi_\xi
\longrightarrow 0.\] 
By sheafification, it becomes
\[
 0 \longrightarrow \mathcal B(0,1,-1)\longrightarrow
\mathcal B\longrightarrow\mathcal B(1,1,\mu_0)\longrightarrow 0, \]
where $\mathcal B$ is $\widetilde{\textrm{Coker}\,A}$ .
We tensorise it with the locally free sheaf $\mathcal B(0,1,-1)$ and we get:
 \[
 0 \longrightarrow \mathcal{O}_Y\longrightarrow
\mathcal B\otimes\mathcal B(0,1,-1)\longrightarrow\mathcal B(1,1,-\mu_0)\longrightarrow 0. \]
Since $\mathcal B((0,1),1,\mu_0)$ is uniquely determined by the existence of the
exact sequence
\[
0 \longrightarrow \mathcal{O}_Y \longrightarrow \mathcal B((0,1),1,\mu_0) \longrightarrow \mathcal B(1,1,-\mu_0)
\longrightarrow 0,   \]
the vector bundle $\mathcal B\otimes\mathcal B(0,1,-1)$ is isomorphic to
$\mathcal B((0,1),1,\mu_0)$. But this means that $\mathcal B$ is isomorphic to $\mathcal B((0,1),1,\mu_0)$.

Consider the map $\varphi_\xi : R(-2)\oplus R(-3) \longrightarrow R(-1)^2$ as in the subsection 1.1, such that the degree of $\textrm{Coker}\,\varphi_\xi$ is -1.\\
From the exact sequence ${\textstyle(3)}$ we see that deg($\Omega^1_R(M_2)$)=$-9$.
Therefore the $R$-module corresponding to $\mathcal B((0,1),1,\mu_0)$ is $(\Omega^1
_R(M_2)\otimes\textrm{Coker}\ \varphi_\xi)^{\vee\vee}\otimes R(2)$. It is the cokernel of the graded map
$A : R(-2)^5 \longrightarrow R(-1)^3 \oplus R \oplus R(-1)$, defined by the matrix $A$ (see above).

(2) The previous statement implies that \\ $\widetilde{\Omega^1_R(M_2)}=\mathcal B((0,1),1,\mu_0)\otimes \widetilde{(\textrm{Coker}\,\varphi_\xi)^\vee}
\otimes \mathcal O_Y(-2)$.\\
Using $\ref{2gen}$ and $\ref{deg1bundle}$ we obtain that $\widetilde{(\textrm{Coker}\,\varphi_\xi)^\vee} =\widetilde{\textrm{Coker}\,\psi_\xi}=\mathcal
B(1,1,\mu_0)$.
Therefore $\widetilde{\Omega^1_R(M_2)}=\mathcal B((0,1),1,\mu_0)\otimes\mathcal B(1,1,\mu_0)\otimes
\mathcal B(-6,1,\mu_0^{-6})= \\\mathcal B((0,1),1,\mu_0)\otimes \mathcal B(-5,1,\mu_0^{-5})=
\mathcal B((-5,-4),1,\mu_0^{-9})$.
   \end{pf}


\subsection{The classification of non--locally free, rank 2, MCM R--modules}

It is known that on a smooth curve any vector bundle of rank $r\geq 2$, say $\ke$, fits in an extension
\[
0\ra\kl\ra\ke\ra\kf\ra 0
\]
where $\kf$ is a vector bundle of rank $r-1$ and $\kl$ is a line  bundle.\\
 Over an isolated curve singularity, any coherent sheaf $\kc$ of rank $r$ has an extension of type
\[
0\ra\kc_1\ra\kc\ra\kc_2\ra 0
\]
where $\kc_1$ and $\kc_2$ are coherent sheaves of  rank 1, respectively $r-1$.
 It is possible that $\kc_1$ and $\kc_2$  are non--locally free but $\kc$ is vector bundle.

The theorem $\ref{mfext}$ describe the extensions of graded MCM modules over a hypersurface ring.
 It gives us an algorithm to compute matrix factorizations of the non--locally free, rank two, indecomposable, graded MCM modules
over the ring $R$. Their minimal number of generators is smaller equal to 6.
The one that are three minimally generated, are isomorphic, up to shiftings, to $\Co \beta_s$, as it was proved in  $\ref{3gen2}$.

\subsubsection{6--generated}

 For each $m\in \Z$, $m \geq 1$, define the matrix:
$$\delta_m=\left(  \begin{matrix}  0 & y_1 & y_2&0&y_3^{m}&-y_3^{m}\\
             y_1 & y_2&0&0&y_3^{m}&-y_3^{m}\\
             y_3 & 0 & -y_1-y_3&0&0&y_3^{m}\\
                    0&0&0&0 & y_1 & y_2\\
                     0&0&0&y_1 & y_2&0 \\
             0&0&0&y_3 & 0& -y_1-y_3\end{matrix}\right).$$

\begin{thm} \label{6sing}
There are countably many isomorphism classes of graded, indecomposable, rank two, non--locally free MCM $R$-modules that
are minimally 6--generated.\\
They are cokernel of graded maps defined by the matrices $$\{\delta_m, \delta^t_m | m\in\Z, m \geq 1\}.$$
\end{thm}

\begin{pf}
Let $M$ be a graded, indecomposable, rank two, non--locally free MCM $R$-modules with $\mu(M)=6$.
Then, up to a shifting, $M$ fits in a graded extension of the type
\begin{equation}\label{6e1}
0\ra \Co \alpha_s \ra M \ra \Co \alpha_\lambda \otimes R(k) \ra 0
\end{equation}
or in one of the type
\begin{equation}\label{6e2}
0\ra \Co \alpha_\lambda \otimes R(k) \ra M \ra \Co \alpha_s \ra 0
\end{equation}
with  $k\in \Z, \lambda=(a:b:1) \in V(f)$.

In the above graded exact sequences we consider, as in the previous section, the graded maps $\alpha_\lambda : R(-2)^3
\longrightarrow R(-1)^3$,  for which $\Co \alpha_\lambda$ has degree 0. So $\Co \alpha_\lambda \otimes R(k)$ has degree
$3k$.

The modules $M$ with  extensions of type ($\ref{6e2}$)  are duals of some modules from the first extension.
Therefore, it is enough to prove the statement for the modules with an extension of type ($\ref{6e1}$).

By Theorem $\ref{mfext}$, $M$ has a matrix factorization $(S,S')$, with $S=\left(\begin{smallmatrix} \alpha_s & D \\0 & \alpha_\lambda \end{smallmatrix}\right)$, where $D$ has homogeneous entries
and it fulfills  $\beta_s \cdot D \cdot \beta_\lambda=0$ mod $(f)$.\\
The corresponding
graded map $S$, is defined as \\$S : R(-2)^3\oplus R(k-2)^3\To R(-1)^3\oplus R(k-1)^3$, so, the matrix $D$ should have homogeneous entries of
degree $1-k$.
If $k\geq 2$, the extension splits, if $k=1$ the module $M$ decomposes. We need, therefore, to consider only the
negative shiftings of $\Co \alpha_\lambda$.\\
Denote the entries of $D$ with $d_1,...,d_9$, so that $D=\left( \begin{smallmatrix}
d_1&d_2&d_3\\
d_4&d_5&d_6\\
d_7&d_8&d_9 \end{smallmatrix}\right)$ and denote $m=1-k$.
The matrix $S$ has the form:

$$S=\left(  \begin{matrix}  0 & y_1 & y_2&d_1&d_2&d_3\\
             y_1 & y_2&0&d_4&d_5&d_6\\
             y_3 & 0 & -y_1-y_3&d_7&d_8&d_9\\
 0&0&0&0 & y_1-ay_3 & y_2-by_3\\
 0&0&0&y_1 & y_2+by_3 & (a^2+a)y_3\\
             0&0&0&y_3 &             0    & -y_1-(a+1)y_3\end{matrix}\right).$$

We make some linear transformations, in order to simplify the matrix $S$.\\
By subtracting the first 3 columns from the last three, (with a corresponding multiplication factor) we can "kill" the
variable $y_1$ in all entries of $D$.
We kill $y_2$ in $d_1$ by subtracting again the third column from the fourth one; the new appeared $y_1$ in $d_7$
disappear if we subtract the line 5 from the third one. \\
We eliminate also $y_3$ in $d_1$, using the line 6 and column 2. So we can consider that $d_1=0$.\\
In the same way, we eliminate $y_2$ and $y_3$ from $d_4$, using the column 2 and line 5, respectively line 6 and column
1.
So we can suppose also that $d_4=0$.\\
Using the first column and the line 5, one can eliminate $y_3$ in $d_7$. Therefore we consider $d_7=a_7y_2^{m}$, with
$a_7$ a constant.
Using the first column and the line 4, one can eliminate $y_3$ also in $d_8$, so $d_8=a_8y_2^{m}, a_8$ constant.\\
 The same, using the second column and the line 4, we eliminate $y_2$ in $d_5$ and we write $d_5=a_5y_3^{m}, a_5\in k$.\\
Notice that, if we denote the third column of $S$ with $c_3$ and first column with $c_1$,
$-c_3+(1-a)c_1$ has on the third position $y_1-ay_3$.  So, if we eliminate $y_2$ from $d_9$ subtracting
the fourth line from the third one, we can make again $a_8y_2^m$ on the position [3,5] subtracting $-c_3+(1-a)c_1$ from the column 5.
We can destroy the new $y_1$ from the position [2,5] with the line 4, but it still remains there a term $g\cdot y_3$, with $g$ a
polynomial of degree $m-1$. If $m=1$, $g$ is a constant, so we do not need to make any other transformations.
If $m \geq 2$, we eliminate the possible $y_2$ from the new $d_5$ using the second column and the line 4.
So, at the end of these transformations, we get $d_9=a_9y_3^m$, $a_9$ constant.

We check now the condition $\beta_s \cdot D \cdot \beta_\lambda=0$  to get more informations on the
entries of $D$. We use the procedure \texttt{condext} that was already defined in the second section.
\begin{verbatim}
 ring R3=0,(y(1..3),d(1..9),a,b),(c,dp(3),dp(9),dp(2));
 ideal i=y(1)^3+y(1)^2*y(3)-y(2)^2*y(3),a3+a2-b2;
 qring S3=std(i); 
 matrix A[3][3]=0,y(1)-a*y(3),     y(2)-b*y(3),
             y(1),y(2)+b*y(3),     (a2+a)*y(3),
             y(3),          0,-y(1)-(a+1)*y(3);
 matrix B=subst(A,a,0,b,0);
 matrix D[3][3]=d(1..9); D[1,1]=0; D[2,1]=0;
 ideal P=condext(B,A,D); P;
P[1]=-y(2)*d(8)-y(3)*d(6)-y(3)*d(7)*a^2-y(3)*d(7)*a+
      y(3)*d(8)*b-y(3)*d(9)*a-y(3)*d(9)
P[2]=y(2)*d(7)*b-y(2)*d(8)*a-y(3)*d(2)*a^2-y(3)*d(2)*a+
     y(3)*d(3)*b+y(3)*d(5)*b-y(3)*d(6)*a
P[3]=y(2)*d(7)*a-y(3)*d(2)*b+y(3)*d(3)*a+y(3)*d(5)*a+
     y(3)*d(7)*a*b-y(3)*d(8)*a^2+y(3)*d(9)*b
P[4]=-y(2)*d(2)+y(2)*d(9)+y(3)*d(2)*b-y(3)*d(3)*a
P[5]=y(2)^2*d(8)+y(2)*y(3)*d(2)*a+y(2)*y(3)*d(2)+y(2)*y(3)*d(6)-
     y(3)^2*d(5)*a^2-y(3)^2*d(5)*a+y(3)^2*d(6)*b
P[6]=y(2)^2*d(7)+y(2)*y(3)*d(3)+y(2)*y(3)*d(5)-y(3)^2*d(5)*b+
     y(3)^2*d(6)*a
P[7]=-y(2)*d(2)+y(2)*d(9)+y(3)*d(2)*b-y(3)*d(3)*a
\end{verbatim}
Since no $y_3$ appear in $d_7$ and $d_8$, from the conditions P[6] and P[5]
we get that $d_7=d_8=0$. Then, P[1]  gives $d_6=-d_9(a+1)$.
\begin{verbatim}
 P=subst(P,d(7),0,d(8),0, d(6),-(a+1)*d(9));
 P=simple(P); P=interred(P); P;
P[1]=d(2)*b-d(3)*a-d(5)*a-d(9)*b
P[2]=d(2)*a^2+d(2)*a-d(3)*b-d(5)*b-d(9)*a^2-d(9)*a
P[3]=y(2)*d(3)+y(2)*d(5)-y(3)*d(5)*b-y(3)*d(9)*a^2-y(3)*d(9)*a
P[4]=y(2)*d(2)-y(2)*d(9)-y(3)*d(2)*b+y(3)*d(3)*a
\end{verbatim}
From the condition P[3], we obtain that $d_3$ is divisible with $y_3^m$.
But $d_3$ is a homogeneous polynomial of degree $m$,  so $d_3=a_3y_3^m$, $a_3$ constant. The condition P[3] becomes
$y_3( -a^2a_9-aa_9-a_5b)+y_2(a_3+a_5)=0.$ So, $a_3 = -a_5$ and $a(a+1)a_9=-a_5b$.
With this information, the ideal P becomes:
\begin{verbatim}
 P=subst(P,d(3),-d(5)); P=simple(P); P;
P[1]=d(2)*b-d(9)*b
P[2]=d(2)*a^2+d(2)*a-d(9)*a^2-d(9)*a
P[3]=d(5)*b+d(9)*a^2+d(9)*a
P[4]=y(2)*d(2)-y(2)*d(9)-y(3)*d(2)*b-y(3)*d(5)*a
\end{verbatim}
From the condition P[4], in a similar way as above, one get $d_2=a_2y_3^m$, with $a_2$ constant and $y_3( aa_5+ba_2)+y_2(a_9-a_2)=0.$
Therefore, $a_2=a_9$ and $ aa_5+ba_9=0$.
We write the new form of $S$:
$$S=\left(  \begin{matrix}  0 & y_1 & y_2&0&a_9y_3^m&-a_5y_3^m\\
             y_1 & y_2&0&0&a_5y_3^m&-(a+1)a_9y_3^m\\
             y_3 & 0 & -y_1-y_3&0&0&a_9y_3^m\\
 0&0&0&0 & y_1-ay_3 & y_2-by_3\\
 0&0&0&y_1 & y_2+by_3 & (a^2+a)y_3\\
             0&0&0&y_3 &             0    & -y_1-(a+1)y_3\end{matrix}\right)$$
with $ aa_5+ba_9=0$.\\
 If $a\neq 0$, $d_5=-d_9b/a$. $\Co S$ do not decomposes, so, $a_9 \neq 0$ and can be chosen to be $a$.
We obtain in this way the matrix
$$S=\left(  \begin{matrix}  0 & y_1 & y_2&0&ay_3^m&by_3^m\\
             y_1 & y_2&0&0&-by_3^m&-(a^2+a)y_3^m\\
             y_3 & 0 & -y_1-y_3&0&0&ay_3^m\\
 0&0&0&0 & y_1-ay_3 & y_2-by_3\\
 0&0&0&y_1 & y_2+by_3 & (a^2+a)y_3\\
 0&0&0&y_3 &   0      & -y_1-(a+1)y_3 \end{matrix}\right).$$
But in this case, $D=y_3^{m-1}(\alpha_s-\alpha_\lambda)$, therefore, $\Co S$ decomposes, still.
Suppose now that $a=b=0$. 
The condition that the module $\Co S$ is not locally free is equivalent to 
Fitt$_2(S)R_{\langle y_1,y_2\rangle} \neq R_{\langle y_1,y_2\rangle}$. (see \cite{Pfister})\\
(Fitt$_2$ is the ideal generated by all 4--minors  of the matrix)\\
To check this, we substitute, in the matrix $S$, the variables  $y_1$ and $y_2$ with 0 and the variable $y_3$ with 1.
The fitting ideal of the new matrix is zero if and only if $\Co S$ is non--locally free module.
We obtain the ideal generated by $a_5^2-a_9^2$, thus,  $a_5=a_9$ or $a_5=-a_9$.\\ 
As before, $a_9 \neq 0$ and can be chosen to be $1$.
With this choice, the matrix $S$ becomes $\delta_m$ or $\delta_m^t$.

We prove now the indecomposability.\\
If $\delta_m$ decomposes, there exist the invertible matrices $U$ and $V$, and there exist $\nu_1, \nu_2\in V(f)$ such that:
$\delta_m \cdot U=V\cdot \left(\begin{smallmatrix}
\alpha_{\nu_1} & 0\\0 &\alpha_{\nu_2}\end{smallmatrix}\right)$.\\
Write $U=  \left(\begin{smallmatrix}U_1 & U_2\\U_3 &U_4\end{smallmatrix}\right)$ and
$V=  \left(\begin{smallmatrix}V_1 & V_2\\V_3 &V_4\end{smallmatrix}\right)$.\\
Let $j-1$ be the degree of the entries of the matrices $U_1$ and $V_1$ and $j\, '-1$, the degree of the entries of $U_2$ 
and
$V_2$.
Then, the entries of $U_3$ and $V_3$ should have degree $j-m$, the one of $U_4$ and $V_4$ should have degree $j\,'-m$.\\
Since $\textrm{det}\ U=\textrm{det}\ V=1$, if $m\geq 2$, $U_3=0$ (and $V_3=0$) or $U_4=0$ (and $V_4=0$).
If $m=1$, since $\alpha_sU_3=V_3\alpha_{\nu_1}$, we get $U_3=V_3=0$ or $\nu_1=s$ and $U_3=V_3=t$Id, $t\neq 0$.
In both cases, (also for $m\geq 2$), we obtain that $\nu_1=\nu_2=s$ and that
$$\left(\begin{matrix}
0&y_3^{m}&-y_3^{m}\\
0&y_3^{m}&-y_3^{m}\\
0&0&y_3^{m}
\end{matrix}\right) =V_2\alpha_s-\alpha_sU_2,$$ that is impossible.
 (in the right hand-side, the entry [1,2] is in the ideal $\langle y_1,y_2 \rangle$,
so, it can not be  $-y_3^{m}$).
Therefore, for any $m$, the matrices $\delta_m$ and $\delta^t_m$ are indecomposable.\\
With a similar proof, one can show that there not exist two invertible matrices $U$ and $V$ such that
$U\delta_m=\delta^t_mV$; more, because of degree reason, it is clear that for two different $m_1$ and
$m_2$, $\delta_{m_1}$ and $\delta_{m_2}$ do not give isomorphic modules. This complete the proof of the theorem.
\end{pf}

Remark: The proof of the indecomposability of the matrices defined in the following two theorems is very similar with the one 
made on the former proof. The computations are simple, but laborious and we decided to skip them.


\subsubsection{5--generated modules}

For all $m\in \Z$, $m\geq 1$, we define the matrices:

$$\alpha^m_{\psi 1}=\left(\begin{matrix}  0 & y_1 & y_2&y_3^m&-y_3^m\\
             y_1 & y_2&0&-y_3^m&y_3^m\\
             y_3 & 0 & -y_1-y_3&y_3^m&0\\
 0&0&0 & y_1^2+y_1y_3&\phantom{-} -y_2y_3\\
             0&0&0&-y_2 &y_1\end{matrix}\right),$$
$$\alpha^m_{\psi 2}=\left(  \begin{matrix}  0 & y_1 & y_2&y_3^m&y_3^m\\
             y_1 & y_2&0&y_3^m&y_3^m\\
             y_3 & 0 & -y_1-y_3&-y_3^m&0\\
            0&0&0 & y_1^2+y_1y_3 & \phantom{-}-y_2y_3\\
             0&0&0&-y_2&y_1\end{matrix}\right),$$
$$\alpha^m_{\varphi 1}=\left(  \begin{matrix}  0 & y_1 & y_2&     y_3^m & y_3^{m+1}\\
                         y_1 & y_2&  0  &     -y_3^m & -y_3^{m+1}\\
                         y_3 & 0  & -y_1-y_3& 0  & y_3^{m+1}\\
                    0&0&0 & y_1 & y_2y_3\\
                    0&0&0& y_2 &\phantom{-} y_1^2+y_1y_3 \end{matrix}\right),$$
$$\alpha^m_{\varphi 2}=\left(  \begin{matrix}  0 & y_1 & y_2&     y_3^m & -y_3^{m+1}\\
                         y_1 & y_2&  0  &     y_3^m & -y_3^{m+1}\\
                         y_3 & 0  & -y_1-y_3& 0  & y_3^{m+1}\\
                    0&0&0 & y_1 & y_2y_3\\
                    0&0&0&\phantom{-} y_2 \phantom{-}& y_1^2+y_1y_3 \end{matrix}\right).$$

For each $\lambda \in V(f)$ we define the matrix $\alpha_{\psi}^{\lambda}$ as follows:\\
If $\lambda=(a:b:1)$ then
$$\alpha_{\psi}^{\lambda}=\left(\begin{matrix}0 & y_1 & y_2&by_3&-ay_3\\
             y_1 & y_2&0&-(a+a^2)y_3&by_3\\
             y_3 & 0 & -y_1-y_3&ay_3&0\\
             0&0&0 & y_1^2+(a+1)y_1y_3+(a^2+a)y_3^2&\  -y_2y_3-by_3^2\\
             0&0&0&-y_2+by_3 &\ \ y_1-ay_3\end{matrix}\right),$$
if $\lambda=(0:1:0)$, 
$$\alpha_{\psi}^{\lambda}=\left(\begin{matrix}0 & y_1 & y_2&0&0\\
             y_1 & y_2&0&y_2&-y_2\\
             y_3 & 0 & -y_1-y_3&0&y_2\\
             0&0&0 & y_1^2&\  -y_2^2\\
             0&0&0&-y_3 &\ \ y_1+y_3\end{matrix}\right).$$
\begin{thm}
The isomorphism classes of graded, indecomposable, rank two, non--locally free MCM $R$-modules that
are minimally 5--generated are given by the matrices:
 $$\{ \alpha^m_{\psi 1}, \alpha^m_{\psi 2}, \alpha^m_{\varphi 1}, \alpha^m_{\varphi 2}, \alpha_{\psi}^\lambda | m\in\Z, m \geq 1,\lambda\in V(f)\}$$
and their transpositions.
\end{thm}
\begin {pf}
Let $M$ be a graded, indecomposable, rank two, non--locally free MCM $R$-modules with $\mu(M)=5$.
Then, up to a shifting, $M$ or $M^\vee$ fits in one of the following graded extensions:

\begin{equation}\label{e1}
0\ra \Co \alpha_s \ra M \ra \Co \psi_\lambda \otimes R(k) \ra 0
\end{equation}
\begin{equation}\label{e2}
0\ra \Co \alpha_\lambda \ra M \ra \Co \psi_s \otimes R(k) \ra 0
\end{equation}
\begin{equation}\label{e3}
0\ra \Co \alpha_s \ra M \ra \Co \varphi_\lambda \otimes R(k) \ra 0
\end{equation}
\begin{equation}\label{e4}
0\ra \Co \alpha_\lambda \ra M \ra \Co \varphi_s \otimes R(k) \ra 0
\end{equation}
with  $k\in \Z, \lambda \in V(f)$.

$\bullet$ Consider first the extension $(\ref{e1})$.

The module $M$ has a matrix factorization $(S,S')$, with $S=\left(\begin{smallmatrix} \alpha_s & D \\0 & \psi_\lambda \end{smallmatrix}\right)$, as in theorem $\ref{mfext}$.\\
Denote the entries of $D$ with $d_1,...,d_6$, so that $D=\left( \begin{smallmatrix}d_1&d_2\\
d_3&d_4\\
d_5&d_6 \end{smallmatrix}\right)$.\\
 Since the corresponding graded map $S$, is
defined as $$S : R(-2)^3\oplus R(k-2)^2\To R(-1)^3\oplus R(k)\oplus R(k-1),$$
the matrix $D$ should have homogeneous entries of degree $m=1-k$.\\
But, if $k\geq 2$, the extension splits, if $k=1$ the module $M$ decomposes. We consider, therefore, only the negative shiftings of $\Co \psi_\lambda$.

As in the previous proof, we make some linear transformations to eliminate the variable
$y_1$ from all $d_i, i=1,...,6$ and the variable $y_2$ from $d_1, d_3$ and $d_5$. We write $d_1=a_1y_3^m$, $d_3=a_3y_3^m$ and $d_5=a_5y_3^m$ with $a_1, a_3, a_5$ constants.

The extension condition,  $\beta_s \cdot D \cdot \varphi_\lambda=0$, implies that $y_2$
do not appear in $d_2$ and $d_4$ and that $d_6=0$. 
\begin{verbatim}
  setring S; //the ring from Theorem 2.6

  matrix psil[2][2]=
      y(1)^2+(a+1)*y(1)*y(3)+(a2+a)*y(3)^2,-(y(2)*y(3)+b*y(3)^2),
                            -(y(2)-b*y(3)),          y(1)-a*y(3);

  matrix phil[2][2]=
      y(1)-a*y(3),                   y(2)*y(3)+b*y(3)^2,
      y(2)-b*y(3), y(1)^2+(a+1)*y(1)*y(3)+(a2+a)*y(3)^2;

  matrix B=subst(A,a,0,b,0); matrix D[3][2]=d(1..6);
  P=condext(B,psil,D); P;
P[1]=y(2)*d(6)-y(3)*d(3)-y(3)*d(5)*a-y(3)*d(5)-y(3)*d(6)*b
P[2]=y(2)*d(2)+y(2)*d(5)-y(3)*d(1)*a-y(3)*d(2)*b
P[3]=y(2)*d(1)-y(2)*d(4)+y(3)*d(3)*a+y(3)*d(4)*b
P[4]=d(1)*b+d(2)*a^2+d(2)*a-d(4)*b+d(5)*a^2+d(5)*a+d(6)*a*b
P[5]=d(1)*a+d(2)*b-d(4)*a+d(5)*b+d(6)*a^2
\end{verbatim}
From P[1], we find $d_6(y_2-by_3)=y_3^{m+1}(a_3+a_5+aa_5)$. Therefore, $d_6=0$ and $a_3=-a_5-aa_5$.
The same, from P[2], we have $d_2=a_2y_3^m$, $a_2$ constant. More, $a_2=-a_5$ and  $aa_1=a_5b$.
 P[3] implies $d_4=a_4y_3^m$ and $a_4=a_1$, $a_4\in K$.
So the matrix $S$ is looking like:
$$S=\left( \begin{matrix}  0 & y_1 & y_2&a_1y_3^m & -a_5y_3^m\\
                         y_1 & y_2 & 0 & -(a_5+aa_5)y_3^m & a_1y_3^m\\  
                         y_3 & 0 & -y_1-y_3 & a_5y_3^m & 0\\
                           0 & 0 & 0 & y_1^2+(a+1)y_1y_3+(a^2+a)y_3^2 & \phantom{-}-y_2y_3-by_3^2\\
                           0 & 0 & 0 &-y_2+by_3 & y_1-ay_3\end{matrix}\right)$$
with $aa_1-ba_5=0$.\\
If $a \neq 0$ then $a_1=ba_5/a$. Since $\Co S$ should not decompose, $a_5$ is nonzero, and it can be chosen to be $a_5=a$.\\
With this choice, for $m\geq 2$, $D$ can be written as $D=-y_3^{m-2}(\alpha_sU_1+U_2\psi_\lambda)$ for 
$U_1=\left(\begin{smallmatrix} -y_1-(a+1)y_3 & 0 \\
                      0 & y_3\\
                      -y_3& 0 \end{smallmatrix}\right)$ and  
$U_2=\left(\begin{smallmatrix} 0 & -y_3 \\
                      1 & 0\\
                      0 & 0 \end{smallmatrix}\right)$ . So, $\Co S$  decomposes.\\
For $m=1$, we obtain the matrix $\alpha_{\psi_\lambda}$, for $\lambda=(a:b:1)$. This matrix is indecomposable, since $D$ can not be written 
as combination of $\lambda_s$ and $\psi_\lambda$.

Consider $a=b=0$. The matrix $S$ provides a non-locally free module if and only if Fitt$_2(S)\cdot R_{\langle y_1,y_2\rangle}\neq R_{\langle y_1,y_2\rangle}$. This means that $a_1^2=a_5^2$. \\Since $a_5$ is nonzero and we can choose it to be 1, we obtain in this way the matrices  $\alpha^m_{\psi 1}$ or $\alpha^m_{\psi 2}.$

If $\lambda=(0:1:0)$, with very similar computations, we obtain only one indecomposable
 extension, $\Co \alpha_\psi^\lambda$.

$\bullet$ In the case of the extension (8), the computations are very similar to the above one and provide the same matrices as the extension (7).

$\bullet$ Consider now the extension ($\ref{e3}$) (the proof and results in the last case are identical with this one). The matrix $S$ has the form 
$$S=\left(  \begin{matrix}  0 & y_1 & y_2&d_1&d_2\\
                           y_1 & y_2& 0  & d_3&d_4\\
                          y_3 & 0 & -y_1-y_3&d_5&d_6\\
 0&0&0& y_1-ay_3 & y_2y_3+by_3^2\\
 0&0&0& y_2-by_3 & y_1^2+(a+1)y_1y_3+(a^2+a)y_3^2\end{matrix}\right).$$
In this case, the corresponding graded map $S$ is defined as \\ $S : R(-2)^3\oplus R(k-2)\oplus R(k-3)\To R(-1)^3\oplus 
R(k-1)^2$.\\
Therefore, on the first column, the matrix $D$ should have homogeneous entries of degree $m=1-k$ and on the second column of degree
$m+1$.
So, if $k\geq 3$, the extension splits, if $k=2$ the module $M$ decomposes.
In the case $k=1$, the first column of $D$ should be 0. The condition $\beta_s \cdot D \cdot \psi_\lambda=0$, implies that also the second column annihilates, so the module $\Co S$ decomposes.\\
Therefore, it is enough to consider  $k\leq 0$, that means, $m\geq 1$.

We make  again some linear transformations to eliminate the variable $y_1$ from all entries of $D$ and the variable $y_2$ from $d_5$. 
Using the fourth line and the second respectively third column, we eliminate  $y_2y_3$ from $d_2$ and $d_6$. 
Subtracting the first column from the fourth one we eliminate also the variable $y_3$ from $d_5$, so we can consider $d_5=0$.
(the new appeared $y_1$ in $d_3$ is killed using the fourth line)

The extension condition $\beta_s \cdot D \cdot \psi_\lambda=0$ gives:
\begin{verbatim}
  D[3,1]=0;
  P=condext(B,phil,D); P;
P[1]=d(4)+d(6)*a+d(6)
P[2]=y(2)*d(4)+y(2)*d(6)*a+y(2)*d(6)-y(3)*d(4)*b-
     y(3)*d(6)*a*b-y(3)*d(6)*b
P[3]=-y(3)*d(1)*b-y(3)*d(3)*a-d(2)*a+d(6)*b
P[4]=y(3)*d(1)*a^2+y(3)*d(1)*a+y(3)*d(3)*b+d(2)*b+d(4)*a
P[5]=y(2)*y(3)*d(3)+y(3)^2*d(3)*b+y(2)*d(2)+y(3)*d(4)*a
P[6]=y(2)*y(3)*d(1)+y(3)^2*d(1)*b-y(2)*d(6)+y(3)*d(2)*a
\end{verbatim}
Therefore $d_4=-d_6(a+1)$. More, from P[5] we see that $y_3 | d_2$ and, since we have eliminated $y_2y_3$ from $d_2$, it must have the form 
$d_2=a_2y_3^{m+1}$, $a_2$ constant.
The same, P[6] implies that $d_6=a_6y_3^{m+1}$, $a_6$ constant. Furthermore, P[6] implies also that $y_3^m | d_1$, that, from degree reasons means, $d_1=a_1y_3^m$. Looking at the coefficients of $y_2y_3^{m+1}$ and $y_3^{m+2}$ we obtain that $a_1=a_6$ and $a_1b+a_2a=0$.
Similarly, from P[5] we get that $d_3=a_3y_3^m$ and $a_3=-a_2$. So, the matrix $S$ has the form:
$$S=\left(  \begin{matrix}  0 & y_1 & y_2     &  a_1y_3^m & a_2y_3^{m+1}\\
                          y_1 & y_2 &  0      & -a_2y_3^m & -a_1(a+1)y_3^{m+1}\\
                          y_3 &  0  & -y_1-y_3&        0  &  a_1y_3^{m+1}\\
                             0&  0  &   0     & y_1-ay_3  & y_2y_3+by_3^2\\
                             0&  0  &   0     & y_2-by_3  & y_1^2+(a+1)y_1y_3+(a^2+a)y_3^2 \end{matrix}\right).$$ 
 If $a \neq 0$ then $a_2=-ba_1/a$.  So $a_1 \neq 0$ and it can be chosen to be $a_1=a$.\\
But then we have $D=-y_3^{m-1}(\alpha_sU_1+U_2\phi_\lambda)$ for 
$U_1=\left(\begin{smallmatrix}  0 & y_1+(a+1)y_3  \\
                           1 &  0\\
                           0 & y_3 \end{smallmatrix}\right)$ and  
$U_2=\left(\begin{smallmatrix} -1 & 0  \\
                           0 & -1\\
                           0 & 0 \end{smallmatrix}\right)$ . Therefore, the module $\Co S$  decomposes.\\
If $a=b=0$, the matrix $S$ provide a non-locally free module if and only if $a_1^2=a_2^2$. Since $a_1$ is nonzero and we can choose it to be 1, we obtain in this way the matrices  
 $\alpha^m_{\varphi 1}$ or $\alpha^m_{\varphi 2}$.\\
With similar computations, one can prove that there are no indecomposable extensions of this type, if 
$\lambda=(0:1:0)$.  \end{pf}


\subsubsection{4--generated modules}

 For all $\lambda=(a:b:1)\in V(f)$, with  $a \neq 0$ we define the matrix:
$$\varphi_{\psi \lambda}=\left(  \begin{matrix}  y_1 & y_2y_3& -by_3&¸ay_3\\
             y_2&\ \ y_1^2+y_1y_3&ay_1+(a^2+a)y_3&-by_3\\
              0&0& \ \ y_1^2+(a+1)y_1y_3 +(a^2+a)y_3^2 & \ \ -y_2y_3-by_3^2\\
             0&0&-y_2+by_3&\ \  y_1-ay_3\end{matrix}\right).$$
If $\lambda=(0:1:0)$ let be:
$$\varphi_{\psi \lambda}=\left(  \begin{matrix}  y_1 & y_2y_3& 0&y_2\\
             y_2&\ \ y_1^2+y_1y_3&y_1&0\\
              0&0& y_1^2 & \phantom{-}-y_2^2\\
             0&0&-y_3& \phantom{-}y_1+y_3 \end{matrix}\right).$$
For all $m\in \Z, m\geq 1$ we define:\\
$$\varphi^m_{\psi 1}=\left(\begin{matrix} y_1 & y_2y_3 & y_3^m & y_3^m\\
             y_2&y_1^2+y_1y_3&y_1y_3^{m-1}+y_3^m & -y_3^m\\
              0&0& y_1^2+y_1y_3 &-y_2y_3\\
             0&0&-y_2& y_1\end{matrix}\right),$$
$$\varphi^m_{\psi 2}=\left(\begin{matrix} y_1 & y_2y_3& -y_3^m &y_3^m\\
             y_2&y_1^2+y_1y_3&y_1y_3^{m-1}+y_3^m &y_3^m\\
              0&0& y_1^2+y_1y_3 &-y_2y_3\\
             0&0&-y_2&y_1\end{matrix}\right),$$
$$\psi^m_{\varphi 1}=\left(\begin{matrix}
 y_1^2+y_1y_3 &- y_2y_3 &y_3^{m+1}& y_1y_3^{m+1}+y_3^{m+2}\\
 -y_2 &y_1& y_3^m&y_3^{m+1}\\
  0 &0& y_1&y_2y_3\\
  0&0&y_2& y_1^2+y_2y_3\end{matrix}\right),$$
$$\psi^m_{\varphi 2}=\left(\begin{matrix} y_1^2+y_1y_3 &- y_2y_3&y_3^{m+1}&-y_1y_3^{m+1}-y_3^{m+2}\\
                               -y_2&y_1&-y_3^m&y_3^{m+1}\\
                           0&0& y_1&y_2y_3\\
                           0&0&y_2&
y_1^2+y_2y_3\end{matrix}\right),$$
$$\varphi^m_{\varphi 1}=\left(\begin{matrix} y_1&y_2y_3&y_3^m&-y_3^{m+1}\\
             y_2& y_1^2+y_1y_3&y_3^m&-y_1y_3^m-y_3^{m+1}\\
              0&0& y_1&y_2y_3\\
             0&0&y_2&y_1^2+y_1y_3\end{matrix}\right),$$
$$\varphi^m_{\varphi2}=\left(\begin{matrix} y_1&y_2y_3&-y_3^m&-y_3^{m+1}\\
             y_2&y_1^2+y_1y_3&y_3^m&y_1y_3^m+y_3^{m+1}\\
             0&0&y_1&y_2y_3\\
             0&0&y_2&y_1^2+y_1y_3\end{matrix}\right).$$
\begin{thm}
The isomorphism classes of graded, indecomposable, rank two, non--locally free MCM $R$-modules that
are minimally 4--generate are given by the matrices $$\{ \varphi^m_{\psi 1}, \varphi^m_{\psi 2}, \psi^m_{\varphi 1}, \psi^m_{\varphi 2},
 \varphi^m_{\varphi 1}, \varphi^m_{\varphi 2}| m\in\Z, m \geq 1\}\cup 
\{ \varphi_{\psi \lambda} | \lambda \in V(f)_{\textrm{reg}}.\}$$ and their transpositions.
\end{thm}
\begin{pf}
Let $M$ be a graded, indecomposable, rank 2, non--locally free MCM module with $\mu(M)=4$.
With a similar argumentation as in previous theorems, up to a shifting, $M$ or $M^\vee$ fits in one of the following graded extensions:

\begin{equation}\label{4e1}
0\ra \Co \varphi_s \ra M \ra \Co \psi_\lambda \otimes R(k) \ra 0
\end{equation}
\begin{equation}\label{4e2}
0\ra \Co \psi_s \ra M \ra \Co \varphi_\lambda \otimes R(k) \ra 0
\end{equation}
\begin{equation}\label{4e3}
0\ra \Co \varphi_s \ra M \ra \Co \varphi_\lambda \otimes R(k) \ra 0
\end{equation}
with  $k\in \Z, \lambda \in V(f)$.

  $\bullet$  Suppose $M$ has an extension of type ($\ref{4e1}$).

Then the module $M$ has a matrix factorization $(S,S')$, with
$S=\left(\begin{smallmatrix} \varphi_s & D \\0 & \psi_\lambda \end{smallmatrix}\right)$ and the matrix  $D$ has homogeneous entries of degree $m=1-k$.\\
So, if $k\geq 2$, the extension splits, if $k=1$ the module $M$ decomposes. We need, therefore, to consider only the
negative shiftings of $\Co \psi_\lambda$.
Denote the entries of $D$ with $d_1,...,d_4$, so that $D=\left( \begin{smallmatrix}d_1&d_2\\
d_3&d_4 \end{smallmatrix}\right)$.\\
Consider first the case $\lambda=(a:b:1)$.
By linear transformations, we can eliminate the variable $y_1$ in $d_1, d_2, d_4$ and $y_2$ in $d_4$.
In case that $m\geq 2$, we can eliminate also $y_2y_3$ in $d_2$ and $y_1^2$ in $d_3$.
We write $d_3=y_1d_5+d_6$ and $d_4=a_4y_3^m$.

The extension condition  $\psi_s \cdot D \cdot \varphi_\lambda=0$ mod $(f)$ :
\begin{verbatim}
 matrix D[2][2]=          d(1),d(2),
                y(1)*d(5)+d(6),d(4);
 phi=subst(phil,a,0,b,0); 
 P=condext(phi,psil,D); P;
P[1]=y(3)*d(5)*a+y(3)*d(5)-d(6)
P[2]=y(2)*d(1)-y(2)*d(4)+y(3)*d(4)*b+y(3)*d(6)*a
P[3]=y(3)*d(5)*b-d(1)*a-d(2)*b+d(4)*a
P[4]=-d(1)*b-d(2)*a^2-d(2)*a+d(4)*b+d(6)*a
P[5]=y(2)*y(3)*d(5)-y(2)*d(2)+y(3)*d(1)*a+y(3)*d(2)*b
P[6]=-y(2)*d(2)*a-y(2)*d(2)+y(2)*d(6)+y(3)*d(4)*a^2+y(3)*d(4)*a+
      y(3)*d(6)*b
\end{verbatim}
We make the substitution $d_3=d_5(y_1+(a+1)y_3)$.(see P[1])
\begin{verbatim}
 P=subst(P,d(6),y(3)*d(5)*a+y(3)*d(5));
 P=interred(P); P;
P[1]=y(3)*d(5)*b-d(1)*a-d(2)*b+d(4)*a
P[2]=y(3)*d(5)*a^2+y(3)*d(5)*a-d(1)*b-d(2)*a^2-d(2)*a+d(4)*b
P[3]=y(2)*d(1)-y(2)*d(4)+y(3)*d(1)*b+y(3)*d(2)*a^2+y(3)*d(2)*a
P[4]=y(2)*y(3)*d(5)-y(2)*d(2)+y(3)*d(1)*a+y(3)*d(2)*b
\end{verbatim}
From P[4] we get that $y_3 | d_2$, and since $y_2y_3$ is eliminated from $d_2$, we can write
$d_2=a_2y_3^m$, $a_2$ constant. P[3] implies then that also $d_1$ has the form $d_1=a_1y_3^m$, $a_1$ constant. Replacing $d_1$ and $d_2$ in P[4]we get  $d_5=a_5y_3^{m-1}$, with $a_5$ constant.
Further more, $a_5=a_2$, $a_4=a_1$ and $aa_1+ba_2=0$.

If $a\neq 0$, then $a_2$ should be nonzero (otherwise the matrix decomposes)
and  can be chosen to be $a$. Then $a_1=b$. If $m=1$ we get the matrix $\varphi_{\psi \lambda}$. But, if $m\geq 2$, the matrix $S$, after some simple linear transformations, decomposes.

In case $a=b=0$, the matrix $S$ corresponds to a non-locally free module iff $a_1^2-a_2^2=0$.
Choosing $a_2=1$, we get the matrices $\varphi^m_{\psi 1}$ and $\varphi^m_{\psi 2}$.

If $\lambda=(0:1:0)$, with similar calculation, we get only one indecomposable extension, $\varphi_{\psi \lambda}$.\\ 

 $\bullet$ Consider now that $M$ has an extension of type ($\ref{4e2}$), that means,
$$
0\ra \Co \psi_s \ra M \ra \Co \varphi_\lambda \otimes R(k) \ra 0.
$$
The module $M$ has a matrix factorization $(S,S')$, with $S=\left(\begin{smallmatrix} \psi_s & D \\0 & \varphi_\lambda \end{smallmatrix}\right)$. Denote the entries of $D$ with $d_1,...,d_4$, so that 
$D=\left( \begin{smallmatrix}d_1&d_2\\                                                                                                d_3&d_4 \end{smallmatrix}\right)$.

In this case, $d_3$ has degree $m=1-k$, $d_1$, $d_4$ have degree  $m+1$  and $d_2$ has degree $m+2$.
So, if $k\geq 3$, the extension splits. If $k=2$ the only nonzero entry of $D$ is on position [1,3]
and has degree 1. If $k=1$ the entry [2,1] should be zero.
The condition of extension, implies that $D=0$, so $S$ decomposes.\\
We consider therefore $m\geq 1$.\\
Let $\lambda=(a:b:1)$ be a point on the nodal curve.\\
By linear transformations over $S$, we eliminate $y_1$ in $d_1, d_3, d_4$, $y_2$ in $d_4$ and $y_1^2$ in $d_2$.
We write $d_2=y_1d_5+d_6$ and $d_4=a_4y_3^m$, $a_4$ constant.
 More, in case that $d_5$ has $y_2y_3$ we eliminate it using the third line.\\
The extension condition  $\varphi_s \cdot D \cdot \psi_\lambda=0$ means:
\begin{verbatim}
 psi=subst(psil,a,0,b,0);
 matrix D[2][2]=d(1),y(1)*d(5)+d(6),
                d(3),          d(4);
 P=condext(psi,phil,D); P;
P[1]=y(3)*d(5)*a+y(3)*d(5)-d(6)
P[2]=-y(3)*d(3)*b+d(1)*a-d(4)*a+d(5)*b
P[3]=y(3)^2*d(3)*a^2+y(3)^2*d(3)*a+y(2)*d(1)-y(2)*d(4)+
     y(3)*d(4)*b
P[4]=y(2)*y(3)*d(3)-y(2)*d(5)+y(3)*d(1)*a+y(3)*d(5)*b
P[5]=y(2)*d(1)-y(2)*d(4)+y(3)*d(1)*b+d(6)*a
P[6]=y(3)*d(5)*b^2-d(6)*a^2
\end{verbatim}
We make the substitution $d_6=d_5(a+1)y_3$. (see P[1])
\begin{verbatim} 
 P=subst(P,d(6),y(3)*d(5)*a+y(3)*d(5));
 P=interred(P); P;
P[1]=y(3)*d(3)*b-d(1)*a+d(4)*a-d(5)*b
P[2]=y(2)*d(1)-y(2)*d(4)+y(3)*d(1)*b+y(3)*d(5)*a^2+y(3)*d(5)*a
P[3]=y(3)^2*d(3)*a^2+y(3)^2*d(3)*a+y(2)*d(1)-y(2)*d(4)+
     y(3)*d(4)*b
P[4]=y(2)*y(3)*d(3)-y(2)*d(5)+y(3)*d(1)*a+y(3)*d(5)*b
\end{verbatim}
From P[4] we see that $y_3 | d_5$, and since $y_2y_3$ is eliminated from $d_5$, we can write
$d_5=a_5y_3^{m+1}$, $a_5$ constant. Since $d_4=a_4y_3^{m+1}$, from P[4],  we get that also
$d_1$ has the form $d_1=a_1y_3^{m+1}$. And, therefore,  from P[2], $d_3=a_3y_3^m$.
Furthermore, $a_4=a_1$, $a_5=a_3$ and $ba_3+aa_1=0$.

If $a\neq 0$, $a_3$ should be nonzero, and we can choose it to be $a$. So
$a_1=a_4=-b$. The matrix $S$ has the form:
$$S=\left(\begin{matrix} y_1^2+y_1y_3 & -y_2y_3& -by_3^{m+1}&ay_1y_3^{m+1}-(a^2+a)y_3^{m+2}\\
                       -y_2&y_1& ay_3^m& -by_3^{m+1}\\
                           0&0& y_1-ay_3 & y_2y_3+by_3^2\\
                           0&0&y_2-by_3&y_1^2+(a+1)y_2y_3+(a^2+a)y_3^2\end{matrix}\right)$$
  and decomposes, after some linear transformations.

If $a=b=0$, the module $\Co S$ is non-locally free if and only if $a_1=a_3$ or $a_1=-a_3$. 
We choose $a_1=1$ and  we get the matrices  $\psi^m_{\varphi 1}$ and $\psi^m_{\varphi 2}$.

If $\lambda=(0:1:0)$, one get no indecomposable extensions of this type.

 $\bullet$ We consider now the last case, when $M$ has an extension of type ($\ref{4e3}$).

Let $(S,S')$ be a matrix factorization of the  module $M$ such that
$S=\left(\begin{smallmatrix} \varphi_s & D \\0 & \varphi_\lambda \end{smallmatrix}\right)$.
As before, denote the entries of $D$ with $d_1,...,d_4$, so that $D=\left( \begin{smallmatrix}d_1&d_2\\
               d_3&d_4 \end{smallmatrix}\right)$.\\ 
 This matrix  has homogeneous entries of degree $m=1-k$ on the first column, and of degree $m+1$ on the second one.
So, if $k\geq 2$, the extension splits. If $k=1$ the first column of $D$ should be zero.
Writing the condition of extension, one get easily that $D=0$, so $S$ decomposes.\\
Therefore, it is enough to consider only the negative shiftings of $\Co \varphi_\lambda$.
Let $\lambda=(a:b:1)$.
By linear transformations, we can eliminate the variable $y_1$ in $d_1, d_2, d_3$ and $y_2$ in $d_1$.
We  eliminate also $y_2y_3$ in $d_2$ and $y_1^2$ in $d_4$ and write $d_4=y_1d_5+d_6$ and $d_1=a_1y_3^m$, with $a_1$ constant.

The extension condition $\psi_s \cdot D \cdot \psi_\lambda=0$ gives :
\begin{verbatim}
 matrix D[2][2]=d(1),d(2),d(3),y(1)*d(5)+d(6);
 P=condext(phi,phil,D); P;
P[1]=y(3)*d(5)*a+y(3)*d(5)-d(6)
P[2]=y(2)*d(1)+y(2)*d(5)-y(3)*d(3)*a-y(3)*d(5)*b
P[3]=y(3)*d(1)*b+y(3)*d(3)*a+y(3)*d(5)*b+d(2)*a
P[4]=y(3)*d(1)*a^2+y(3)*d(1)*a+y(3)*d(3)*b+d(2)*b+d(6)*a
P[5]=y(2)*y(3)*d(3)-y(3)^2*d(1)*a^2-y(3)^2*d(1)*a+
     y(2)*d(2)-y(3)*d(2)*b
\end{verbatim}
We make the substitution $d_6=y_3(a+1)d_5$.
\begin{verbatim}
 P=subst(P,d(6),y(3)*d(5)*a+y(3)*d(5)); P=interred(P); P;
P[1]=y(3)*d(1)*b+y(3)*d(3)*a+y(3)*d(5)*b+d(2)*a
P[2]=y(3)*d(1)*a^2+y(3)*d(1)*a+y(3)*d(3)*b+y(3)*d(5)*a^2+
     y(3)*d(5)*a+ d(2)*b
P[3]=y(2)*d(1)+y(2)*d(5)-y(3)*d(3)*a-y(3)*d(5)*b
P[4]=y(2)*y(3)*d(3)-y(3)^2*d(1)*a^2-y(3)^2*d(1)*a+y(2)*d(2)-
     y(3)*d(2)*b
\end{verbatim}
From P[4] we see that $y_3 | d_2$, and since $y_2y_3$ is eliminated from $d_2$, we can write
$d_2=a_2y_3^{m+1}$, $a_2$ constant. Furthermore, since $d_1=a_1y_3^m$, we get that also
$d_3$ has the form $d_3=a_3y_3^m$, with $a_3$ constant. P[1] implies that $d_5=a_5y_3^m$, $a_5$ constant.
More, $a_5=-a_1$, $a_2=-a_3$ and $aa_3-ba_1=0$.

If $a\neq 0$, as before, $\Co S$ decomposes.\\
If $a=b=0$, the module $\Co S$ is non-locally free iff $a_1=a_3$ or $a_1=-a_3$. If we choose $a_3=1$, the matrix $S$ becomes one of $\varphi^m_{\varphi 1}$ or $\varphi^m_{\varphi 2}$.\\
If $\lambda=(0:1:0)$, one get no indecomposable extensions of this type.
\end{pf}

\end{document}